\documentclass[12pt]{article}
\usepackage{amsmath}
\usepackage{amsfonts}
\usepackage{amssymb}
\usepackage{amssymb}
\usepackage{graphicx}
\usepackage{subcaption}
\usepackage{xcolor}
\usepackage{geometry}
\usepackage{array}
\usepackage{booktabs}
\usepackage{zymacros}%
\begin{document}

\title{FBSDE based neural network algorithms for high-dimensional quasilinear
parabolic PDEs}
\author{Wenzhong Zhang\thanks{Department of Mathematics, Southern Methodist
University, Dallas, TX 75275.}
\and Wei Cai\thanks{Corresponding author, Department of Mathematics, Southern
Methodist University, Dallas, TX 75275(\texttt{cai@smu.edu}).} }
\maketitle

\begin{abstract}
In this paper, we propose forward and backward stochastic differential
equations (FBSDEs) based deep neural network (DNN) learning algorithms for the
solution of high dimensional quasilinear parabolic partial differential
equations (PDEs), which are related to the FBSDEs by the Pardoux-Peng theory. The algorithms rely on a learning process by minimizing
the pathwise difference between two discrete stochastic processes,
defined by the time discretization of the FBSDEs and the DNN representation of
the PDE solutions, respectively. The proposed algorithms are shown to generate DNN solutions for a  100-dimensional Black--Scholes--Barenblatt equation,
accurate in a finite region in the solution space, and has a convergence rate similar to that of the Euler--Maruyama discretization used for the FBSDEs. As a result, a Richardson extrapolation technique over time discretizations can be used to enhance the accuracy of the DNN solutions.
For time oscillatory solutions, a multiscale DNN is shown to improve the performance of the FBSDE DNN for high frequencies.

\end{abstract}

\makeatletter \makeatother


\section{Introduction}

The relationship between stochastic processes and the solution of partial
differential equations represents one of the high achievements of probability
theory in potential theory research \cite{doob}, represented by the
celebrated Feynman--Kac formula in linear parabolic and elliptic PDEs as a
result of the Kolmogorov backward equation for the generator of the stochastic
process for the former \cite{pavliotis} and the Dynkin formula for the latter
\cite{oksendal}. The recent work by Pardoux--Peng \cite{pardoux-peng} has
extended the concept of the classic linear Feynman--Kac formula to a nonlinear
version, which connects the solution of a quasilinear parabolic PDE to a
coupled pair of forward and backward stochastic processes. This extraordinary
development has made much impact in the mathematical finance in option pricing \cite{karoui}.

Meanwhile, in the field of scientific computing, this connection between SDEs and quasilinear PDEs has inspired
new approaches of solving high dimensional parabolic partial differential
equations (PDEs), which are ubiquitous in material sciences such as the
Allen--Cahn equation for phase transition, and quantum mechanics such as
the Schrodinger equation as well as the Black--Scholes equation for option pricing
and the Hamilton-Jacobi-Bellman equation for optimal control. For PDEs in high dimensions, the main challenge of the traditional
numerical methods, such as finite element, finite difference and spectral
methods, is the curse of dimensionality, namely, the number of the
unknowns in the discretized systems for the PDEs grows exponentially in terms
of the dimension of the problem. Recently, machine learning approaches based
on the deep neural network have taken advantage of the Pardoux--Peng's theory
for forward and backward stochastic differential equations (FBSDEs) and PDEs.
The solution to the PDEs can be learned by sampling the paths of involved
stochastic processes, which are discretized in time by the classic
Euler--Maruyama scheme \cite{kloeden}. The first such an attempt was done in
\cite{han}, where neural network was used as an approximator to
the gradient of the PDEs solutions, while the PDE's solution follows the
dynamics of the FBSDEs, and the learning was carried out by imposing the
terminal condition provided by the parabolic PDEs. Another approach
\cite{raissi} is to approximate the PDE's solution itself by a deep neural
network, which also provides the gradient of the solution as required by the
FBSDEs, the learning is then carried out by minimizing the difference between
the solution given by the discretized SDEs and that given by the DNN at all
discretization time stations. In this paper, improved learning schemes will
be proposed based on a similar approach in \cite{raissi}, but with clearer
mathematical reasoning for the learning processes, to ensure
the numerical methods' mathematical consistency and improved convergence for the PDEs' solutions.

The rest of the paper is organized as follows. In Section 2, we will review
the Pardoux--Peng's theory, which establishes the relation between FBSDEs and
quasilinear parabolic PDEs, with an emphasis on the relation between the
classic Feynman--Kac formula and the nonlinear version represented by the
Pardoux--Peng theory. Section 3 will first review the algorithms
proposed in \cite{han} and \cite{raissi}, and then two new improved methods will be proposed.
Section 4 will present numerical results of the new schemes for solving a
100-dimensional Black--Scholes--Barenblatt equation. Enhanced numerical accuracy by
Richardson extrapolations and multi-scale DNNs for PDEs with oscillatory solutions in time will also be discussed.
 Finally, a conclusion
will be given in Section 5.

\section{Pardoux--Peng theory on FBSDEs and quasilinear parabolic PDEs}

In this paper, we consider {the scalar solution} $u(t,x),$ $t\in\lbrack0,T]$,
$x\in\mathbb{R}^{d}$ for the following $d$-dimensional{ parabolic} PDE
\begin{equation}
\partial_{t} u+\frac{1}{2}\mathrm{Tr}[\sigma\sigma^{T}\nabla\nabla u]+\mu\cdot\nabla
u=\phi, \label{eq-paraPDE}%
\end{equation}
with a {terminal condition}
\begin{equation}
u(T,x)=g(x),
\end{equation}
where $\sigma=\sigma(t,x,u)$, $\phi=\phi(t,x,u,\nabla u)$, $\mu=\mu
(t,x,u,\nabla u)$ are functions with ranges in  with dimensions $d\times d$, $1$ and $d$,
respectively. We are interested in finding the {initial value} $u(0, x_0)$
given $x_0\in\mathbb{R}^{d}$. Therefore, in some sense our problem is similar
to a time reverse problem for a time reversed version of \ref{eq-paraPDE} with an initial data at $t=0$.

Following Pardoux--Peng in \cite{pardoux-peng}, under certain regularity conditions, the forward-backward SDE reformulation
gives a nonlinear implicit
Feynman--Kac formula for the solution of the parabolic PDE (\ref{eq-paraPDE}).
The FBSDEs are proposed as follows. Let $W_{t}=(W_{t}^{1},\cdots,W_{t}^{d})$
where each $W_{t}^{j}$ is a standard Brownian motion. Let $\{\mathcal{F}%
_{t}:0\leq t\leq T\}$ be its natural filtration on the time interval $[0,T]$.
Then, we have the equations of stochastic processes $X_{t}$, $Y_{t}$ and
$Z_{t}$ in $d$, $1$ and $d$ dimensions that are adaptive to the filtration
$\{\mathcal{F}_{t}:0\leq t\leq T\}$, respectively,
\begin{align}%
\begin{split}
dX_{t}  &  =\mu(t,X_{t},Y_{t},Z_{t})dt+\sigma(t,X_{t},Y_{t})dW_{t},\\
X_{0}  &  = x_0,
\end{split}
\label{fsde}\\%
\begin{split}
dY_{t}  &  =\phi(t,X_{t},Y_{t},Z_{t})dt+Z_{t}^{T}\sigma(t,X_{t},Y_{t}%
)dW_{t},\\
Y_{T}  &  =g(X_{T}).
\end{split}
\label{bsde}%
\end{align}
If $\mu$ and $\sigma$ do not explicitly depend on $Y_{t}$ or $Z_{t}$, the
FBSDEs are \emph{decoupled}.

We can easily show that processes defined by
\begin{equation}
Y_{t}=u(t,X_{t}),\quad Z_{t}=\nabla u(t,X_{t}) \label{SDE-PDE}%
\end{equation}
in fact satisfy the above equations (\ref{fsde}) and (\ref{bsde}).

By using the
Ito's formula \cite{oksendal} and the forward SDE of $X_{t}$, we have
\begin{equation}%
\begin{split}
dY_{t}  &  =du(t,X_{t})=\partial_{t}udt+\nabla u\cdot dX_{t}+\frac{1}{2}%
\sum_{i=1}^{d}\sum_{j=1}^{d}\partial_{ij}u { } d[X^{i},X^{j}]_{t}\\
&  =\partial_{t}udt+\nabla u\cdot\left(  \mu dt+\sigma\cdot dW_{t}\right)
+\frac{1}{2}\mathrm{Tr}[{\sigma}{\sigma}^{T}\nabla\nabla u]dt\\
&  =\left(  \partial_{t}u+\nabla u\cdot{\mu}+\frac{1}{2}\mathrm{Tr}[{\sigma
}{\sigma}^{T}\nabla\nabla u]\right)  dt+Z_{t}^{T}{\sigma}d{W}_{t},
\end{split}
\label{pdeProof}%
\end{equation}
which gives the PDE (\ref{eq-paraPDE}) by comparing (\ref{pdeProof}) with the
backward SDE (\ref{bsde}) for $Y_{t}$.

\noindent

The determination of the third stochastic process $Z_{t}$ from the two SDEs in
(\ref{fsde}) and (\ref{bsde}) makes use of the martingale representation
theory \cite{karatzas}. Consider the following special case of the backward SDE (\ref{bsde})
as an example:
\begin{equation}
Y_{t}+\int_{t}^{T}f(s,X_{s})ds+\int_{t}^{T}Z_{s}\cdot dW_{s}=g(X_{T}%
),\quad0\leq t\leq T, \label{eq-BSDE-v1}%
\end{equation}
i.e. $\mu(t,x,u,\nabla u)=f(t,x)$, and $\sigma(t,x,u)=I_{d\times d}$ is the
identity matrix. By taking the conditional expectation with respect to
$\mathcal{F}_{t}$, we have
\begin{equation}
{Y_{t}=\mathbb{E}\left[  \left.  Y_{t}\right\vert \mathcal{F}_{t}\right]
}=\mathbb{E}\left[  \left.  g(X_{T})-\int_{t}^{T}f(s,X_{s})ds\right\vert
\mathcal{F}_{t}\right]  ,\quad0\leq t\leq T.
\end{equation}
Next, we define the following martingale
\begin{equation}
L_{t}=\mathbb{E}\left[  \left.  g(X_{T})-\int_{0}^{T}f(s,X_{s})ds\right\vert
\mathcal{F}_{t}\right]  ,\quad0\leq t\leq T,
\end{equation}
where $L_{0}=Y_{0}$. By the{ martingale representation theorem \cite{karatzas}%
}, there exists a process ${Z}_{t}^{\star}$ such that
\begin{equation}
L_{t}=Y_{0}+\int_{0}^{t}Z_{s}^{\star}\cdot dW_{s},\quad0\leq t\leq T.
\label{MartinRep}%
\end{equation}
The stochastic process $Z_{t}^{\star}$ is unique in the sense that
\begin{equation}
\int_{0}^{T}\left\Vert {Z}_{t}^{\star}-{Z}_{t}^{\ast}\right\Vert
^{2}dt=0,\quad\text{a.s.}%
\end{equation}
if ${Z}_{t}^{\ast}$ satisfies the same condition (\ref{MartinRep})\ as
${Z}_{t}^{\star}$ \cite{karatzas}.

Meanwhile, we can show that $Z_{t}%
=Z_{t}^{\star}$ solves the backward SDE (\ref{eq-BSDE-v1}),
\begin{align*}
{}  &  Y_{t}+\int_{t}^{T}f(s,X_{s})ds+\int_{t}^{T}Z_{s}^{\star}\cdot
dW_{s}-g(X_{T})\\
={}  &  {\mathbb{E}\left[  \left.  g(X_{T})-\int_{t}^{T}f(s,X_{s}%
)ds\right\vert \mathcal{F}_{t}\right]  }
  +\left(  \int_{0}^{T}{-\int_{0}^{t}}\right)  f(s,X_{s})ds+(L_{T}%
-L_{t})-g(X_{T})\\
={}  &  \int_{0}^{T}f(s,X_{s})ds+L_{T}-g(X_{T})\\
={}  &  L_{T}-\mathbb{E}\left[  \left.  g(X_{T})-\int_{0}^{T}f(s)ds\right\vert
\mathcal{F}_{T}\right] \\
={}  &  0.
\end{align*}


Connection with the classic Feynman--Kac formula is interpreted as follows. If
in the parabolic PDE (\ref{eq-paraPDE}), $\phi$ has a linear dependence on $u$,
i.e.
\begin{equation}
\phi(t,x,u,\nabla u)=c(t,x)u(t,x)+f(t,x),
\end{equation}
then, the backward SDE (\ref{bsde}) has an explicit solution
\begin{equation}%
\begin{split}
Y_{t}={}  &  e^{-\int_{t}^{T}c(s,X_{s})ds}g(X_{T})-\int_{t}^{T}e^{-\int
_{t}^{s}c(\tau,X_{\tau})d\tau}f(s,X_{s})ds 
 {}-\int_{t}^{T}e^{-\int_{t}^{s}c(\tau,X_{\tau})d\tau}Z_{s}^{T}%
\sigma(s,X_{s},Y_{s})dW_{s}.
\end{split}
\label{ExplicitSol}%
\end{equation}
By taking the conditional expectation on both sides, we arrive at
\begin{equation}
Y_{t}=\mathbb{E}\left[  \left.  e^{-\int_{t}^{T}c(s,X_{s})ds}g(X_{T})-\int
_{t}^{T}e^{-\int_{t}^{s}c(\tau,X_{\tau})d\tau}f(s,X_{s})ds\right\vert
\mathcal{F}_{t}\right]  .
\end{equation}
For $(t,x)\in\lbrack0,T]\times\mathbb{R}^{d}$, using $X_{t}=x$ as the initial
condition of the forward SDE (\ref{fsde}) on the time interval $[t, T]$
\emph{instead of} $X_{0} = x_0$, the traditional Feynman--Kac formula
\cite{oksendal} is recovered,
\begin{equation}
u(t,x)=\mathbb{E}\left[  \left.  e^{-\int_{t}^{T}c(s,X_{s})ds}g(X_{T}%
)-\int_{t}^{T}e^{-\int_{t}^{s}c(\tau,X_{\tau})d\tau}f(s,X_{s})ds\right\vert
X_{t}=x\right]  .
\end{equation}

For a general parabolic equation with a nonlinear function $\phi(s,x,u,\nabla
u)$, we have
\[
Y_{t}=\mathbb{E}\left[  \left.  g(X_{T})-\int_{t}^{T}\phi(s,X_{s},Y_{s}%
,Z_{s})ds\right\vert \mathcal{F}_{t}\right]  ,
\]
and for given $(t,x)\in\lbrack0,T]\times\mathbb{R}^{d}$, the following
nonlinear equation for $u(t,x)$ is obtained
\begin{equation}
u(t,x)=\mathbb{E}\left[  \left.  g(X_{T})-\int_{t}^{T}\phi(s,X_{s}%
,u(s,X_{s}),\nabla u(s,X_{s}))ds\right\vert X_{t}=x\right]  .
\end{equation}

\section{FBSDE based neural network algorithms for quasilinear parabolic PDEs}

The learning of the solution will be based on the sample paths of the
FBSDEs, which are linked to the PDE solution in (\ref{SDE-PDE}). Paths of
the FBSDEs will be produced by a {time discretization algorithm with }samples
of the Brownian motion $W_{t}$.

Let $0=t_{0}<\cdots<t_{N}=T$ be a uniform
partition of $[0,T]$. On each interval $[t_{n},t_{n+1}]$, define time and
Brownian motion increments as
\begin{equation}
\Delta t_{n}=t_{n+1}-t_{n},\quad\Delta W_{n}=W_{t_{n+1}}-W_{t_{n}}.
\end{equation}
Denoting $X_{t_{n}}$, $Y_{t_{n}}$ and $Z_{t_{n}}$ by $X_{n}$, $Y_{n}$ and
$Z_{n},$ respectively,  and applying the{ Euler--Maruyama scheme} to the FBSDEs
(\ref{fsde}) and (\ref{bsde}), respectively, we have
\begin{align}
X_{n+1}  &  \approx X_{n}+\mu(t_{n},X_{n},Y_{n},Z_{n})\Delta t_{n}%
+\sigma(t_{n},X_{n},Y_{n})\Delta W_{n},\label{eq-Xn}\\
Y_{n+1}  &  \approx Y_{n}+\phi(t_{n},X_{n},Y_{n},Z_{n})\Delta t_{n}+Z_{n}%
^{T}\sigma(t_{n},X_{n},Y_{n})\Delta W_{n}. \label{eq-Yn}%
\end{align}

Due to the relationship with the parabolic PDE, the solution to the parabolic PDE provides an alternative representation
for $Y_{n+1}$ and $Z_{n+1}$,
\begin{align}
Y_{n+1}  &  =u(t_{n+1},X_{n+1}),\label{eq-Yn-nn}\\
Z_{n+1}  &  =\nabla u(t_{n+1},X_{n+1}). \label{eq-Zn}%
\end{align}

In this paper, fully connected networks of $L$ hidden layers will be used,
which are given in the following form,
\begin{equation}
f_{\vtheta}(\vx)=\vW^{[L-1]}\sigma\circ(\cdots(\mW^{[1]}\sigma\circ
(\mW^{[0]}(\vx)+\vb^{[0]})+\vb^{[1]})\cdots)+\vb^{[L-1]}, \label{dnn}%
\end{equation}
where $W^{[1]}, \cdots, W^{[L-1]}$ and $b^{[1]}, \cdots, b^{[L-1]}$ are the weight
matrices and bias unknowns,respectively, denoted collectively by $\theta$, , to be
optimized via the training, $\sigma(x)$ is the activation function and $\circ$ is the
application of the activation function $\sigma$ applied to a vector quantity component-wisely.

\subsection{Existing FBSDE based neural network algorithms}

\subsubsection{\noindent{Deep BSDE \cite{han}}}

The Deep BSDE trains a network to approximate the random value $Y_N$ at time $t=T$, where $X_{0}=x_0$ is the
input. $Y_{0},Z_{0}$ are trainable variables and $Y_0$ is the targeted quantity of the algorithm.
$W_{n},X_{n},0\leq n\leq N$ can be obtained similarly as before.
The algorithm can be organized as follows.
\begin{enumerate}
\item The initial value $X_{0}=x_0$ is given. Trainable variables $Y_{0}$ and $Z_{0}$ are randomly
initialized.

\item On each time interval $[t_{n},t_{n+1}]$, use the Euler--Maruyama scheme
to calculate $X_{n+1}$ and $Y_{n+1}$ as in (\ref{eq-Xn}) and (\ref{eq-Yn}). Then, train a
fully connected feedforward network
\begin{equation}
f_{\theta}^{(n+1)}(\cdot) \approx \nabla u(t_{n+1}, \cdot)
\end{equation}
where $f^{(n+1)}_{\theta}(\cdot)$ is a fully connected neural network of $H$ hidden layers of the form given in (\ref{dnn}).
Activation functions including ReLU, Tanh, Sigmoid, etc. can been used.

\item Connect all quantities (subnetworks $f^{(n)}_{\theta}(\cdot)$, etc) at $\{ t_n \}$ to form a network that
outputs $Y_{N}$, which is expected to be an approximation of $u(t_{N},X_{N})$.

\item The loss function is then defined by a Monte Carlo approximation of
\begin{equation}
\mathbb{E}\left\|  Y_{N}-g(X_{N})\right\|  ^{2}.
\end{equation}

\end{enumerate}

The Deep BSDE has been shown to give convergent numerical results for various
high dimensional parabolic equations \cite{han} and a posteriori estimate
suggests strong convergence of half order \cite{han2020}.

\bigskip\noindent\textbf{Remark 1}. The Deep BSDE method from \cite{han}
trains the network for the specific initial data $X_{0}=x_0$ and yield only an
approximation to the PDE solution $Y_{0}=u(0,x_0)$. Therefore, once the
desired initial data is changed, a new training may have to be carried out.
Also, the total size of $N$ individual sub-networks used to approximate
$Z_{n}=\nabla u(t_{n},X_{n})$, $n=1,\cdots,N-1$ will grow linearly in terms of
$\ $time discretization steps $N$, resulting in large amount of training
parameter if higher accuracy of the PDE solution is desired.

\bigskip

\subsubsection{FBSNNs \cite{raissi} (\textbf{Scheme 1})}\label{section-FBSNNs}

The FBSNNs trains a network $u_{\theta}(t,x)$ that directly approximates the
solution to the PDE (\ref{eq-paraPDE}) in some region in the $(t, x)$ space.
The network has a fixed size of number of hidden layers and neurons per layer.
The algorithm can be organized as follows.
\begin{enumerate}
\item The initial value $X_{0}=x_0$ is given. Evaluate $Y_{0}$ and $Z_{0}$ using the
network
\begin{equation}
Y_{0}=u_{\theta}(t_{0},X_{0}),\quad Z_{0}=\nabla u_{\theta}(t_{0},X_{0}).
\end{equation}
The gradient above is calculated by an automatic differentiation.

\item On each time interval $[t_{n},t_{n+1}]$, use the Euler--Maruyama scheme
(\ref{eq-Xn}) to calculate $X_{n+1}$, and use the network for $Y_{n+1}$ and
$Z_{n+1}$, i.e.
\begin{align}
\label{scheme1Y1}%
\begin{split}
X_{n+1}  & = X_{n}+\mu(t_{n},X_{n},Y_{n},Z_{n})\Delta t_{n}+\sigma(t_{n}%
,X_{n},Y_{n})\Delta W_{n},\\
Y_{n+1}  & = u_{\theta}(t_{n+1},X_{n+1}),\\
Z_{n+1}  & = \nabla u_{\theta}(t_{n+1},X_{n+1}).
\end{split}
\end{align}
On the other hand, calculate a reference value $Y_{n+1}^{\star}$ using the Euler--Maruyama
scheme (\ref{eq-Yn})
\begin{equation}
Y_{n+1}^{\star}=Y_{n}+\phi(t_{n},X_{n},Y_{n},Z_{n})\Delta t_{n}+Z_{n}%
^{T}\sigma(t_{n},X_{n},Y_{n})\Delta W_{n}. \label{scheme1Y2}%
\end{equation}

\item The loss function is taken as a Monte Carlo approximation of
\begin{equation}
\mathbb{E}\left[  \sum_{n=1}^{N}{\left\Vert Y_{n}-Y_{n}^{\star}\right\Vert
^{2}}+\left\Vert Y_{N}-g(X_{N})\right\Vert ^{2}+\left\Vert Z_{N}-\nabla
g(X_{N})\right\Vert ^{2}\right]  .\label{scheme1loss}%
\end{equation}
In this paper, we will name the above numerical method Scheme 1.
In order to compare the training results using different values
of $N$, the loss function for Scheme 1 is modified as
\begin{align}\label{loss-1}
L_{1}[u_{\theta};x_0]  & =\frac{1}{M}\left[  \sum_{\omega}\frac{1}{N}%
\sum_{n=1}^{N}{\left\Vert Y_{n}-Y_{n}^{\star}\right\Vert ^{2}}+\beta
_{1}\left\Vert Y_{N}-g(X_{N})\right\Vert ^{2}
+\beta_{2}\left\Vert Z_{N}-\nabla g(X_{N})\right\Vert ^{2}\right]
\end{align}
where $M$ serves as the batch size of the training and $\omega$ denotes any
instance of sampling of the discretized Brownian motion $W_{n},0\leq n\leq N-1$,
and $\beta_{1}$, $\beta_{2}$ are the penalty parameters for the terminal
conditions. The averaging factor $1/N$ is introduced for consistency consideration as the reduction of
the loss function as $N$ increases, when applied to the exact solution,
is expected.
\end{enumerate}

\bigskip

\noindent\textbf{Remark 2}. The FBSNNs algorithm proposed in \cite{raissi}
relies on a loss function involving the difference between sequences
$\{Y_{n}\}$ and $\{Y_{n}^{\ast}\}$, which carry the information inside the
time interval $(0, T)$. While the discrete stochastic process $\{Y_{n}\}$ can be
expected to approach a continuous stochastic process as defined in the
backward SDE (\ref{bsde}), the question whether the discrete sequence of
random variables $\{Y_{n}^{\ast}\}$ will converge to the same stochastic process
is not clear. As a result, the rate and extent for the
difference between $\{Y_{n}\}$ and $\{Y_{n}^{\ast}\}$, thus the loss function,
approaching to zero is not certain. Our numerical test will provide some
evidence for this concern.


\noindent

\subsection{Improved FBSDE based deep neural network algorithms for quasilinear parabolic PDEs}

In this section, we propose improved algorithms for the FBSDEs based deep
neural networks similar to the approach in \cite{raissi}, but are
mathematically consistent in the definition of the loss function and the
discretization of both forward and backward SDEs related to the PDE solutions.
Specifically, the loss will be made of the difference of {\it two discrete
stochastic processes}, which will approach the same process given by the
backward SDEs if the overall scheme converges.

\subsubsection{FBSDE based algorithms - Scheme 2}

Based on the Remark 2 from Section \ref{section-FBSNNs}, we would like to design a new scheme
whose loss function is expected to show the strong convergence rate of the
Euler--Maruyama scheme for the discretization of the FBSDEs. A key factor will
be to make the loss function as the pathwise differences between two stochastic
processes, which will converge to the same continuous adapted diffusion
process if the time discretization of FBSDEs and DNN approximations converge.

\bigskip

\noindent\textbf{Scheme 2}. Train a DNN $u_{\theta}(t,x)$ to approximate the
solution $u(t,x)$ of the parabolic PDE (\ref{eq-paraPDE}).

\begin{enumerate}
\item Given $X_{0}=x_0$ and let $Y_{0}=u_{\theta}(t_{0},X_{0})$, $Z_{0}=\nabla
u_{\theta}(t_{0},X_{0})$.

\item On each time interval $[t_{n},t_{n+1}]$, calculate $X_{n+1}$ and
$Y_{n+1}$ using the Euler--Maruyama scheme (\ref{eq-Xn}) and (\ref{eq-Yn}),
respectively, and calculate $Z_{n+1}$ using the network, i.e.
\begin{equation}%
\begin{split}
X_{n+1}  & =X_{n}+\mu(t_{n},X_{n},Y_{n},Z_{n})\Delta t_{n}+\sigma(t_{n}%
,X_{n},Y_{n})\Delta W_{n},\\
{Y_{n+1}}  & =Y_{n}+\phi(t_{n},X_{n},Y_{n},Z_{n})\Delta t_{n}+Z_{n}^{T}%
\sigma(t_{n},X_{n},Y_{n})\Delta W_{n},\\
Z_{n+1}  & = \nabla u_{\theta}(t_{n+1},X_{n+1}).
\end{split}
\label{scheme2Y1}%
\end{equation}
Next, calculate a reference quantity by the DNN representation of the PDE
solution,
\begin{equation}
{Y_{n+1}^{\star}=u_{\theta}(t_{n+1},X_{n+1}).} \label{scheme2Y2}%
\end{equation}

\item For a batch size $M$ with $\omega$ denoting any of the $M$ sample paths,
the loss function is given as

\begin{align}\label{scheme2loss}
L_{2}[u_{\theta};x_0]  & =\frac{1}{M}\sum_{\omega}\left[  \frac{1}{N}%
\sum_{n=1}^{N}\left\Vert Y_{n}-Y_{n}^{\star}\right\Vert ^{2}+\beta
_{1}\left\Vert {Y_{N}^{\star}}-g(X_{N})\right\Vert ^{2}
+\beta_{2}\left\Vert Z_{N}-\nabla g(X_{N})\right\Vert ^{2}\right],
\end{align}
where $\beta_{1}$, $\beta_{2}$ are the penalty parameters of the terminal condition.
\end{enumerate}

The reference quantity $Y_{N}^{\star}$ is used in the terminal term in the loss
function $L_{2}[u_{\theta};x_0]$, because here it is a straightforward output
of the neural network $u_{\theta}$.

\subsubsection{FBSDE based algorithms - Scheme 3}

In the Scheme 2 above, the discrete process (\ref{scheme2Y2}) is defined
through the composite function using the DNN representation of the PDE
solution $u_{\theta}(t,x)$. An alternative way is given below where both
discrete processes are obtained from an Euler--Maruyama discretization of the SDEs.

\bigskip

\noindent\textbf{Scheme 3: }Train a DNN $u_{\theta}(t,x)$ to approximate the
solution $u(t,x)$ of the parabolic PDE (\ref{eq-paraPDE}).

\begin{enumerate}
\item Given the initial values $X_{0}^{(1)}=X_{0}^{(2)}=x_0$ and we compute
\begin{equation}
Y_{0}^{(1)}=Y_{0}^{(2)}=u_{\theta}(t_{0}, x_0),\quad Z_{0}^{(1)}=Z_{0}%
^{(2)}=\nabla u_{\theta}(t_{0}, x_0)
\end{equation}
from the network $u_{\theta}(t,x)$.

\item On each time interval $[t_{n},t_{n+1}]$, calculate $X_{n+1}%
^{(1)},Y_{n+1}^{(1)}$ and $Z_{n+1}^{(1)}$ as in (\ref{scheme1Y1})\ of Scheme
1, then $X_{n+1}^{(2)},Y_{n+1}^{(2)}$ and $Z_{n+1}^{(2)}$ as in
(\ref{scheme2Y1}) of Scheme 2, i.e.
\begin{align}%
\begin{split}
X_{n+1}^{(1)}  &  {}={}X_{n}^{(1)}+\mu(t_{n},X_{n}^{(1)},Y_{n}^{(1)}%
,Z_{n}^{(1)})\Delta t_{n}+\sigma(t_{n},X_{n}^{(1)},Y_{n}^{(1)})\Delta W_{n},\\
Y_{n+1}^{(1)}{}  &  {}=u_{\theta}(t_{n+1},X_{n+1}^{(1)}),\\
Z_{n+1}^{(1)}  &  {}={}\nabla u_{\theta}(t_{n+1},X_{n+1}^{(1)}),
\end{split}
\label{scheme3Y1}\\%
\begin{split}
X_{n+1}^{(2)}  &  {}={}X_{n}^{(2)}+\mu(t_{n},X_{n}^{(2)},Y_{n}^{(2)}%
,Z_{n}^{(2)})\Delta t_{n}+\sigma(t_{n},X_{n}^{(1)},Y_{n}^{(1)})\Delta W_{n},\\
Y_{n+1}^{(2)}{}  &  {}{}=Y_{n}^{(2)}+\phi(t_{n},X_{n}^{(2)},Y_{n}^{(2)}%
,Z_{n}^{(2)})\Delta t_{n}+(Z_{n}^{(2)})^{T}\sigma(t_{n},X_{n}^{(2)}%
,Y_{n}^{(2)})\Delta W_{n},\\
Z_{n+1}^{(2)}  &  {}={}\nabla u_{\theta}(t_{n+1},X_{n+1}^{(2)}).
\end{split}
\label{Scheme3Y2}%
\end{align}

\item For a batch size $M$ with $\omega$ denoting any of the $M$ sample paths,
the loss function is defined by
\begin{align}\label{loss-3}
L_{3}[u_{\theta};x_0]  & =\frac{1}{M}\sum_{\omega}\left[  \frac{1}{N}%
\sum_{n=1}^{N}{\left\Vert Y_{n}^{(1)}-Y_{n}^{(2)}\right\Vert ^{2}}+\beta
_{1}\left\Vert Y_{N}^{(1)}-g(X_{N}^{(1)})\right\Vert ^{2}
+\beta_{2}\left\Vert Z_{N}^{(1)}-\nabla g(X_{N}^{(1)})\right\Vert^{2}\right],
\end{align}
where $\beta_{1}$, $\beta_{2}$ are the penalty parameters of the terminal condition.
\end{enumerate}

\section{Numerical results}
In this section, we will carry out several tests on Scheme 1 from \cite{raissi} and the new Scheme 2
and Scheme 3, for a 100-dimensional Black--Scholes--Barenblatt equation and its variants.
\subsection{100-dimensional Black--Scholes--Barenblatt equation}
Consider the following 100-dimensional Black--Scholes--Barenblatt (BSB) equation from \cite{raissi} as the model problem: for $t \in[0, T]$ and $x\in\mathbb{R}^{d}$, the
scalar function $u(t,x)$ satisfies
\begin{align}
\label{eq-BSB}%
\begin{split}
u_{t} + \frac{1}{2}\mathrm{Tr}\left[  \sigma^{2}\mathrm{diag}(x x^T)\nabla\nabla u\right] &  = r(u-\nabla u\cdot x),\\
u(T,x)  &  = \lVert x\rVert^{2}.
\end{split}
\end{align}
The PDE is linked to the FBSDEs
\begin{align}
\begin{split}
dX_{t}  &  =\sigma\mathrm{diag}(X_{t})dW_{t},\\
X_{0}  &  =x_0,\\
dY_{t}  &  =r(Y_{t}-Z_{t}\cdot X_{t})dt+\sigma Z_{t}^{T}\mathrm{diag}%
(X_{t})dW_{t},\\
Y_{T}  &  =g(X_{T}),
\end{split}
\end{align}
where $g(x)=\lVert x\rVert^{2}$, and $x_0\in\mathbb{R}^{d}$ is the position
where we like to get the initial value $u(0, x_0)$. The exact solution to the PDE
(\ref{eq-BSB}) is given in a closed form by
\begin{equation}
u(t,x)=e^{(r+\sigma^{2})(T-t)}\left\Vert x\right\Vert ^{2},
\end{equation}
so that we can test the accuracy of the DNN schemes.
Parameters are given by $d=100$, $T=1.0$, $\sigma=0.4$, $r=0.05$
and
\begin{equation}
x_0=(1,0.5,1,0.5,\cdots,1,0.5).
\end{equation}

We use a 6-layer fully connected feedforward neural network for
$u_{\theta}(t, x)$ with 5 hidden layers, each having 256 neurons. The
activation function is the sine function as suggested by \cite{raissi}. We
train the network with the Adam optimizer with descending learning rates 1e-3,
1e-4, 1e-5, 1e-6 and 1e-7,
each for 10000 steps.
The batch size is $M=100$.

In the loss functions (\ref{loss-1}), (\ref{scheme2loss}) and (\ref{loss-3}),
the penalty parameters are chosen as $\beta_{1} = \beta_{2} = 0.02$.

Illustration of the training results in the high-dimensional space is provided along the sample paths.
When the training is finished, we randomly generate 1000 sample paths for verification of the accuracy, with a finer time discretization with time steps $\Delta t_{n} = 1/1000$.
For each (discretized) sample path $\omega$ and for $0 \le n \le1000$, the relative error
of this model problem at $(t_{n}, X_{n}(\omega))$ (or at $(t_{n}, X_{n}^{(2)}(\omega))$ when
using Scheme 3) is defined by
\begin{equation}
e_{n}(\omega)=\frac{\lvert u_{\theta}(t_{n}, X_{n}(\omega))-u(t_{n},X_{n}(\omega))\rvert}{\lvert
u(t_{n},X_{n}(\omega))\rvert}.
\end{equation}
The mean and the standard deviation (SD) of each $e_{n}$ can also be calculated.

\subsubsection{Scheme 1 from \cite{raissi}}

\begin{figure}[h]
\centering
\begin{subfigure}{.45\textwidth}
\includegraphics[width=\linewidth]{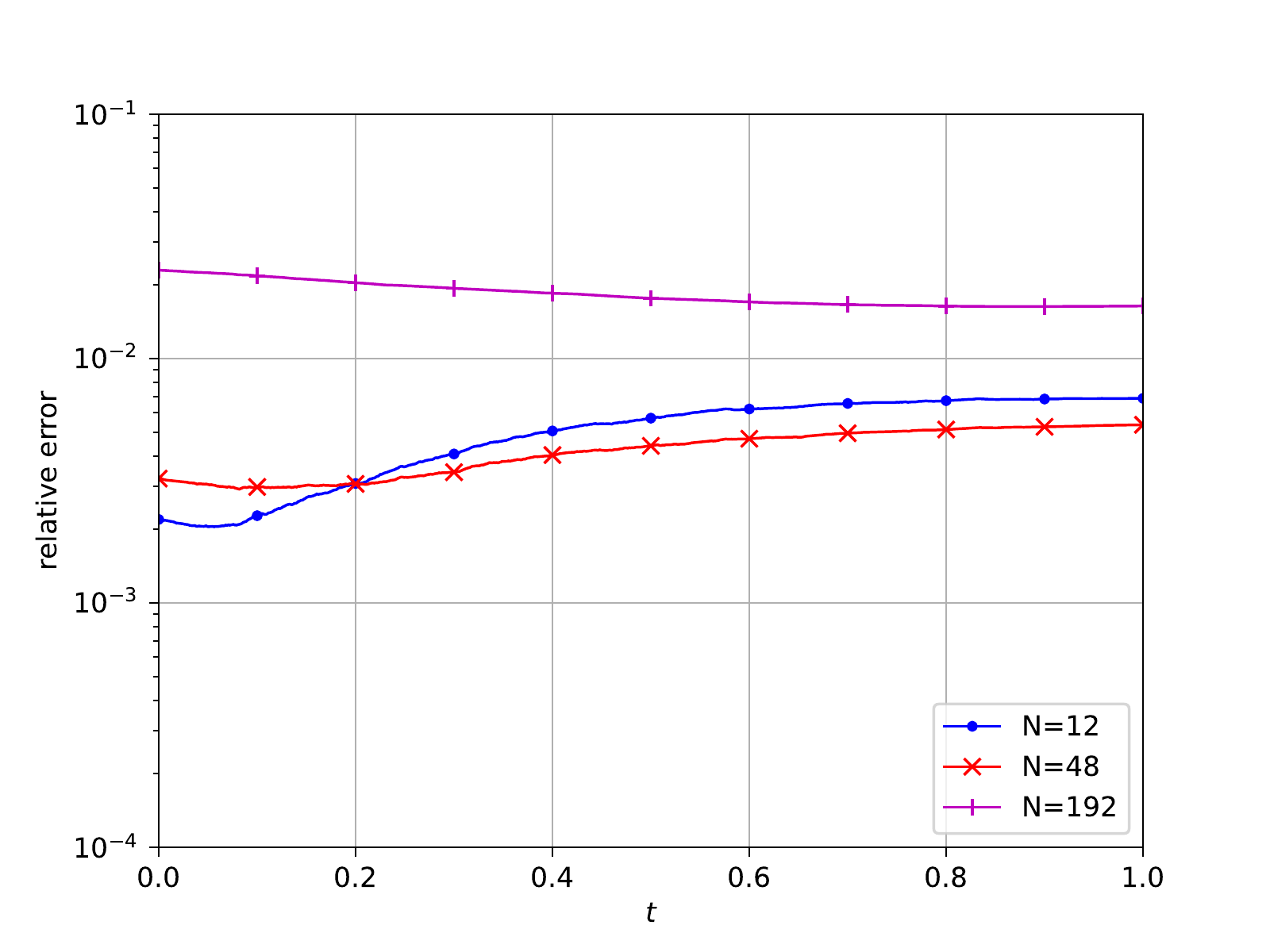}
\caption{Error Mean}
\end{subfigure}
\begin{subfigure}{.45\textwidth}
\includegraphics[width=\linewidth]{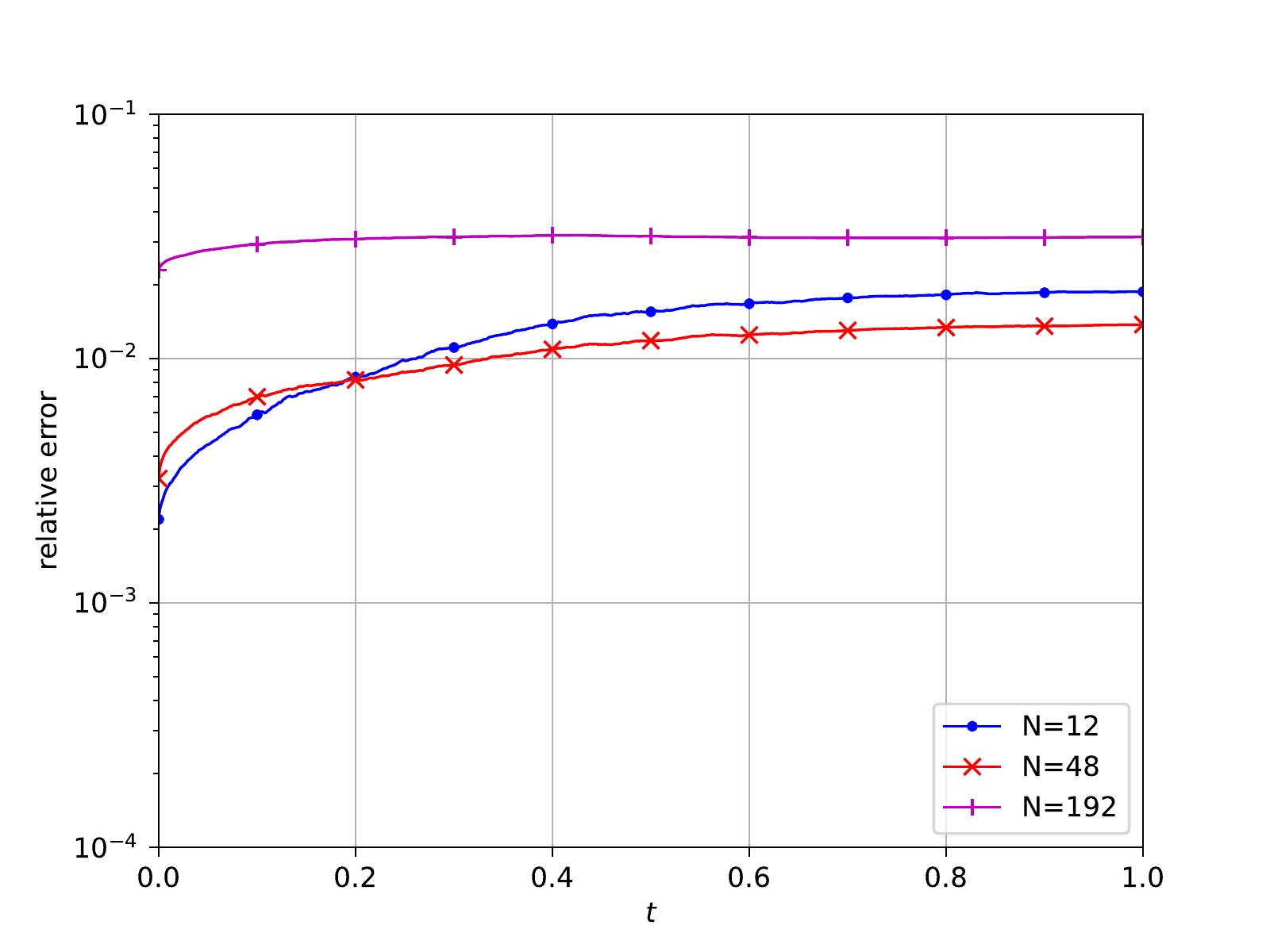}
\caption{Error Mean plus two SDs}
\end{subfigure}
\caption{(Non-convergence) Relative error of Scheme 1 for $N=12$ (middle), $48$ (bottom) and
$192$ (top).}%
\label{error_scheme1}%
\end{figure}

Fig. \ref{error_scheme1} shows the relative error of Scheme 1 for $N=12$, $48$
and $192$, where the mean error and the mean error plus two standard
deviations of the error are presented. We can see the reduction of the errors from
$N=12$ to $N=48$, however, the error increases from $N=48$ to $N=192$. This
degeneracy in accuracy is an indication that as the time discretization is
refined, the two quantities in the definition of loss function
(\ref{scheme1loss}) do not approach the same continuous stochastic process. In
fact, as it is defined by (\ref{scheme1Y2}), $\{Y_{n}^{\star}\}$ may not
converge to a continuous stochastic process at all.

\subsubsection{Scheme 2 and Scheme 3}\label{numerical-S2S3}

\begin{figure}[h]
\centering
\begin{subfigure}{.45\textwidth}
\includegraphics[width=\linewidth]{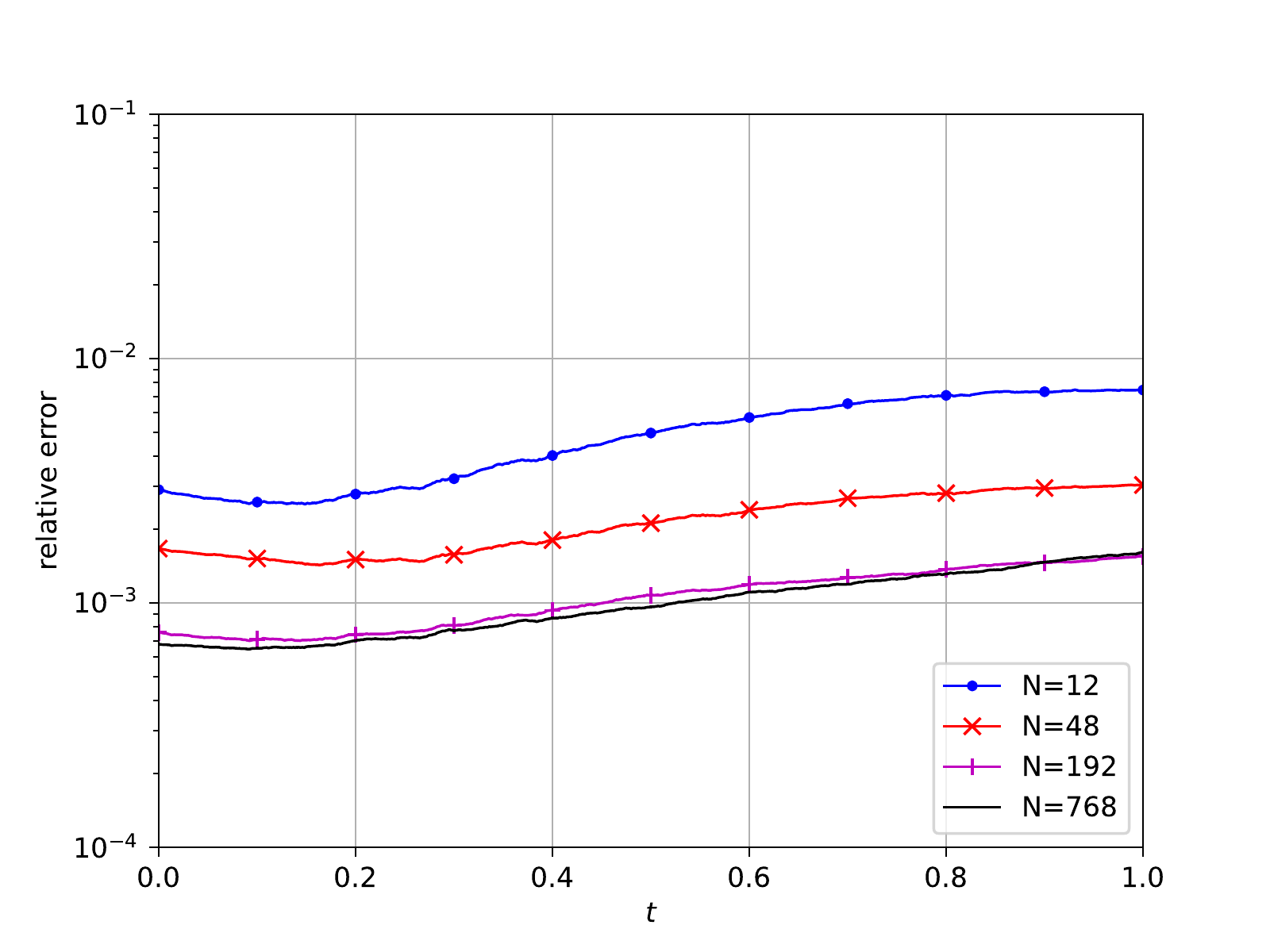}
\caption{Error Mean}
\end{subfigure}
\begin{subfigure}{.45\textwidth}
\includegraphics[width=\linewidth]{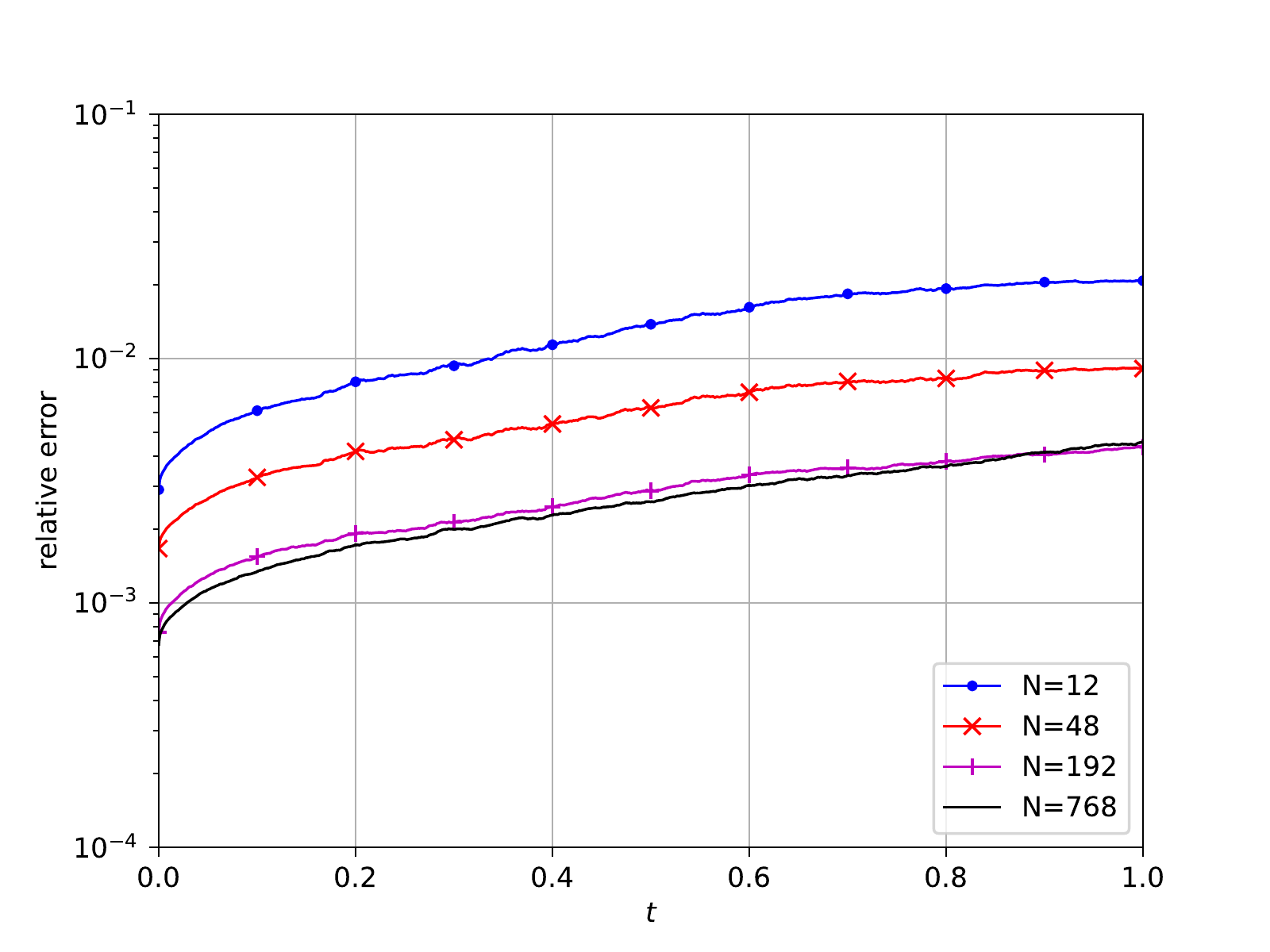}
\caption{Error Mean plus two SDs}
\end{subfigure}
\caption{Relative error of Scheme 2 for $N = 12, 48, 192$ and $768$.}%
\label{error_scheme2}%

\centering
\begin{subfigure}{.45\textwidth}
\includegraphics[width=\linewidth]{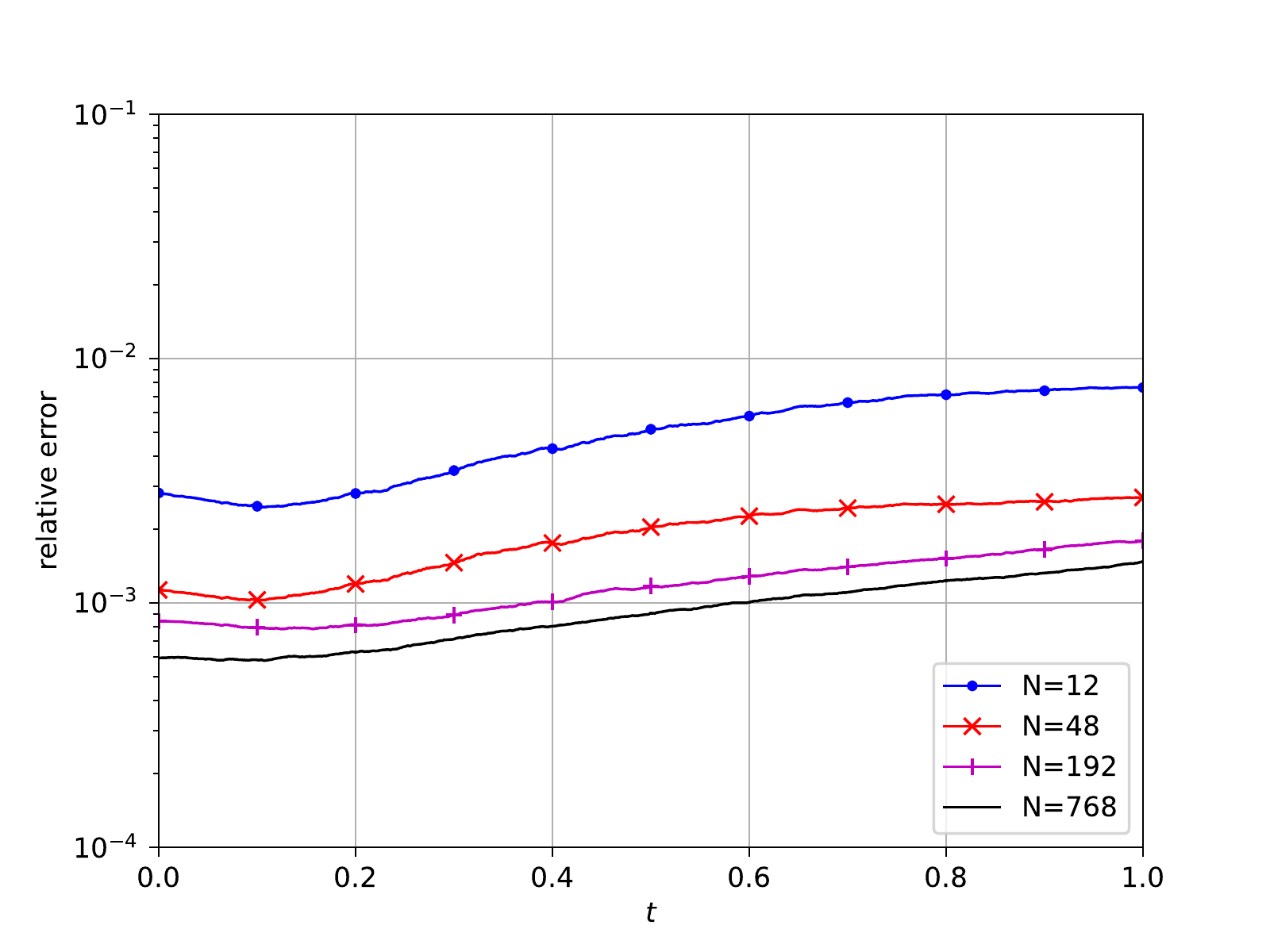}
\caption{Error Mean}
\end{subfigure}
\begin{subfigure}{.45\textwidth}
\includegraphics[width=\linewidth]{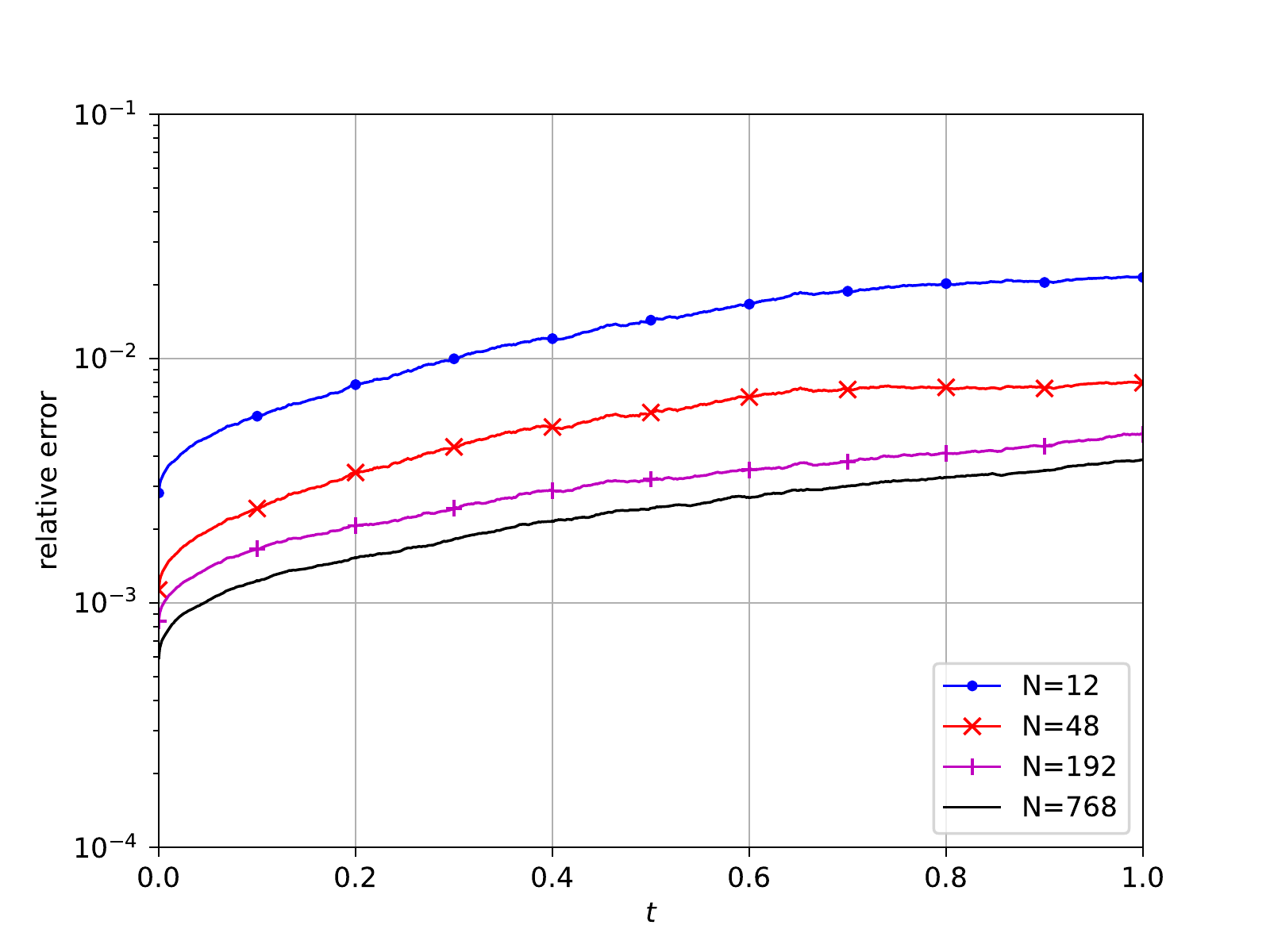}
\caption{Error Mean plus two SDs}
\end{subfigure}
\caption{Relative error of Scheme 3 for $N = 12, 48, 192$ and $768$.}%
\label{error_scheme3}%
\end{figure}

Fig. \ref{error_scheme2} and Fig. \ref{error_scheme3} show the mean error and
mean error plus two standard derivations of the error for Scheme 2 and Scheme
3 for $N=12$, $N=48$, $N=192$ and ${N=768}$, respectively. 

Both the results in Fig. \ref{error_scheme2} and Fig. \ref{error_scheme3} show the
convergence of the new Scheme 2 and Scheme 3, respectively, in contrast to the degeneracy of
the accuracy of Scheme 1 when the time discretization is refined. For both new
schemes, we can see improvement of the accuracy from $N=48$ to $N=192$ is
close to the one from $N=12$ to $N=48$, but the improvement of $N=768$ over
$N=192$ is a little less. This indicates the network training might dominate
the error compared to the time discretization error. In fact, the terminal parts of the loss function failed to halve in the $N=768$ cases compared to $N=192$.

\begin{figure}[hbt!]
\centering
\begin{subfigure}{.45\textwidth}
\includegraphics[width=\linewidth]{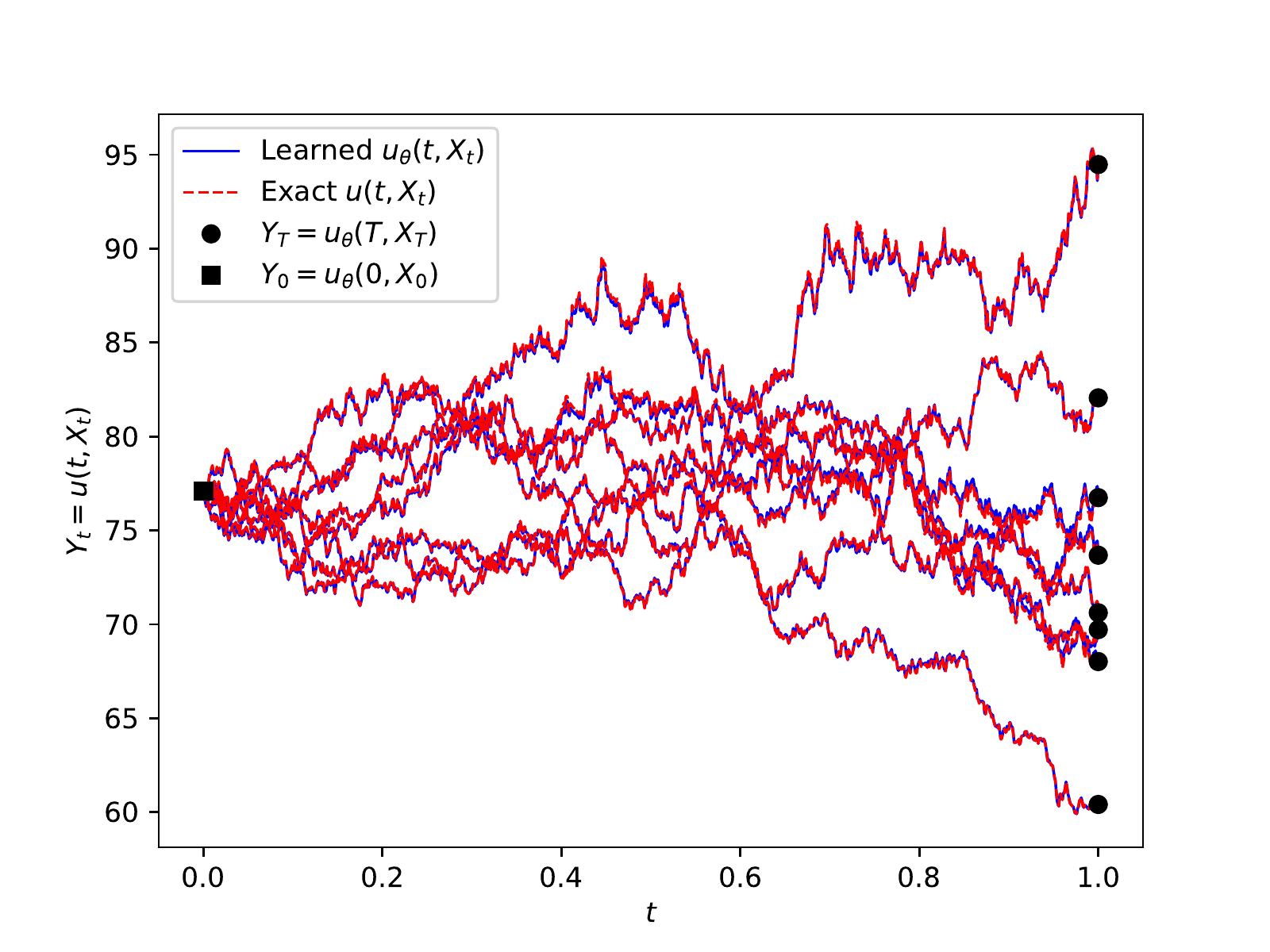}
\caption{Prediction by Scheme 2}
\end{subfigure}
\begin{subfigure}{.45\textwidth}
\includegraphics[width=\linewidth]{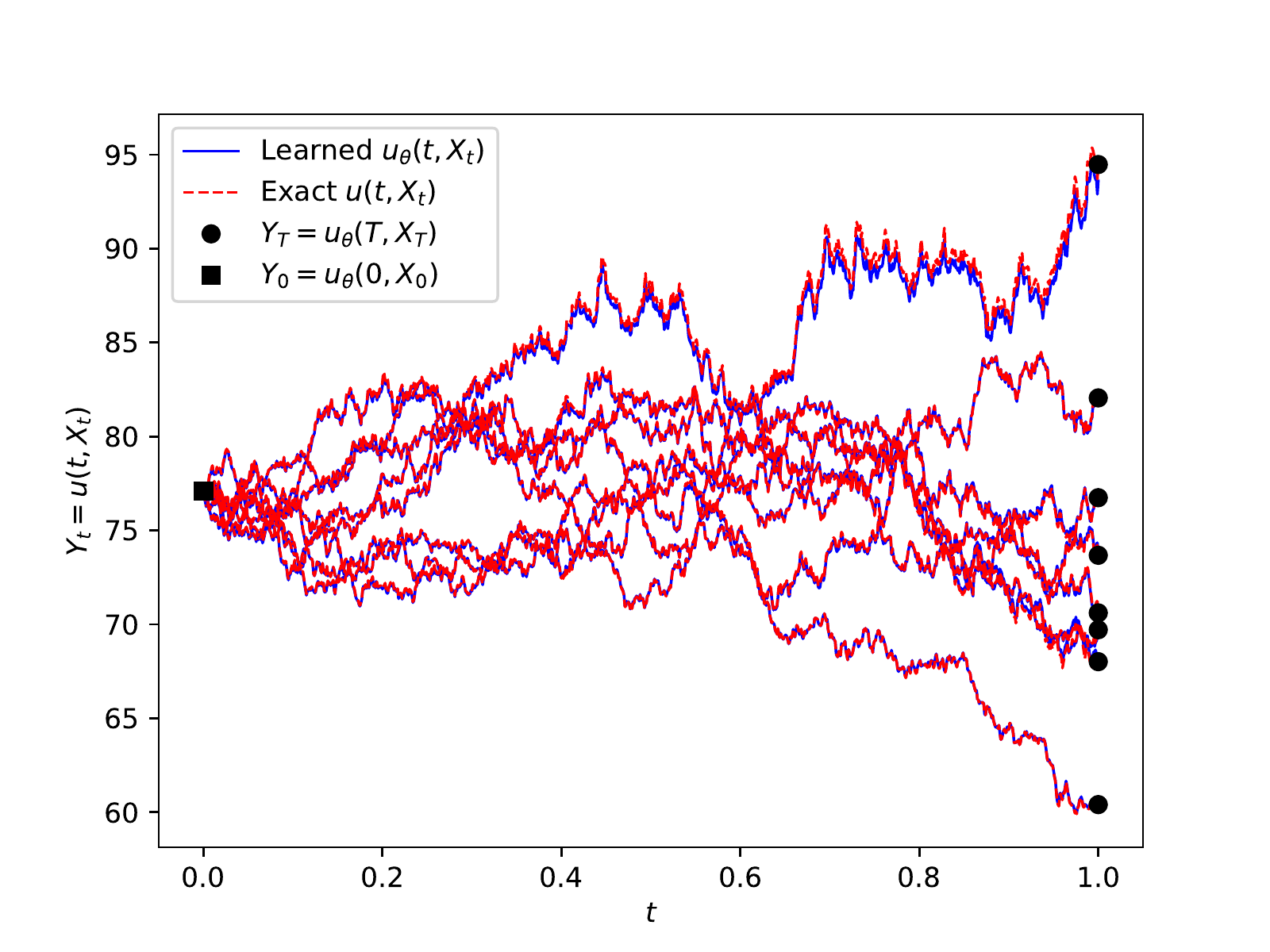}
\caption{Prediction by Scheme 3}
\end{subfigure}
\hfill\begin{subfigure}{.45\textwidth}
\includegraphics[width=\linewidth]{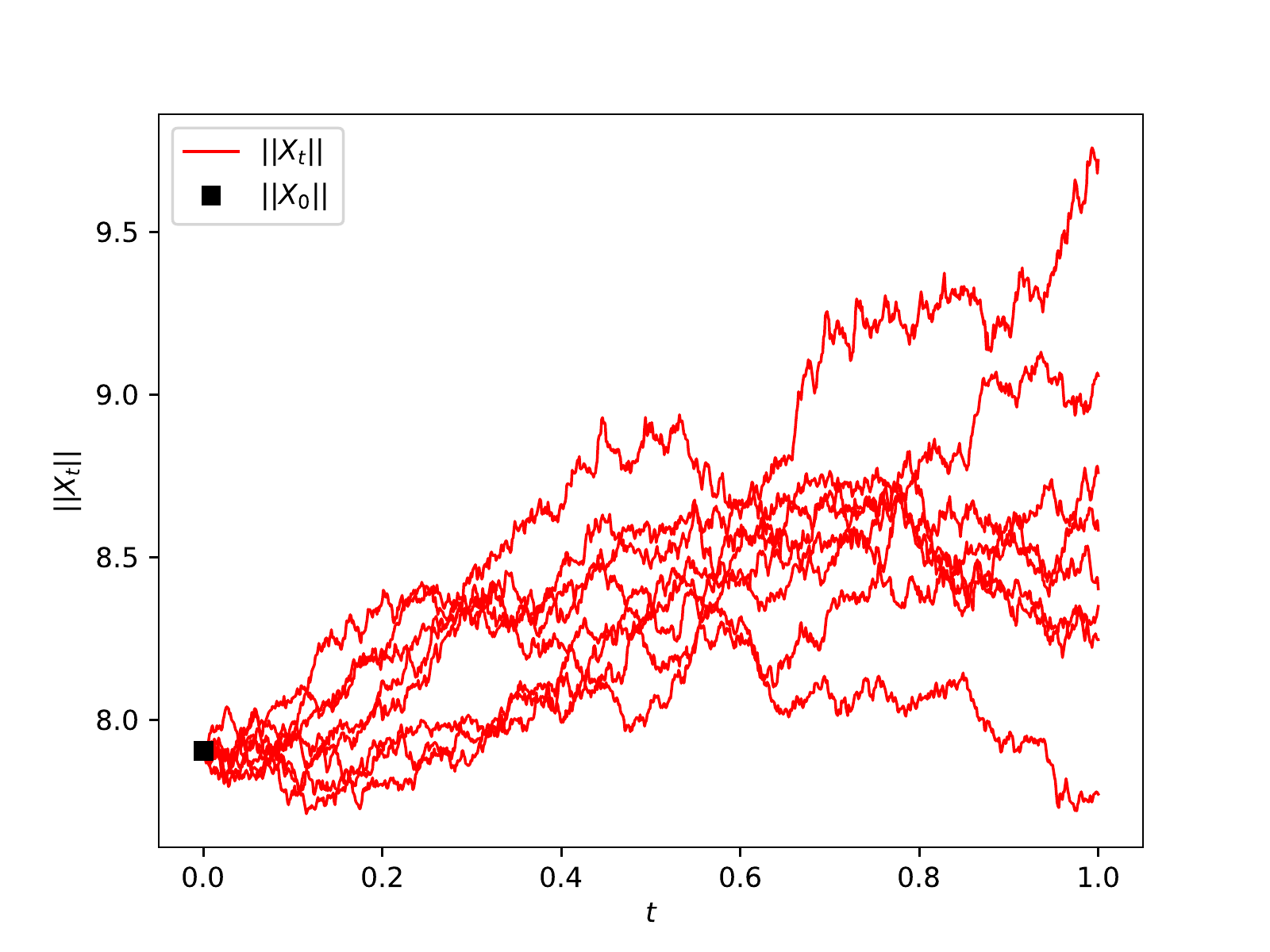}
\caption{$\lVert X_t \rVert$ along 8 sample paths}
\end{subfigure}
\begin{subfigure}{.45\textwidth}
\includegraphics[width=\linewidth]{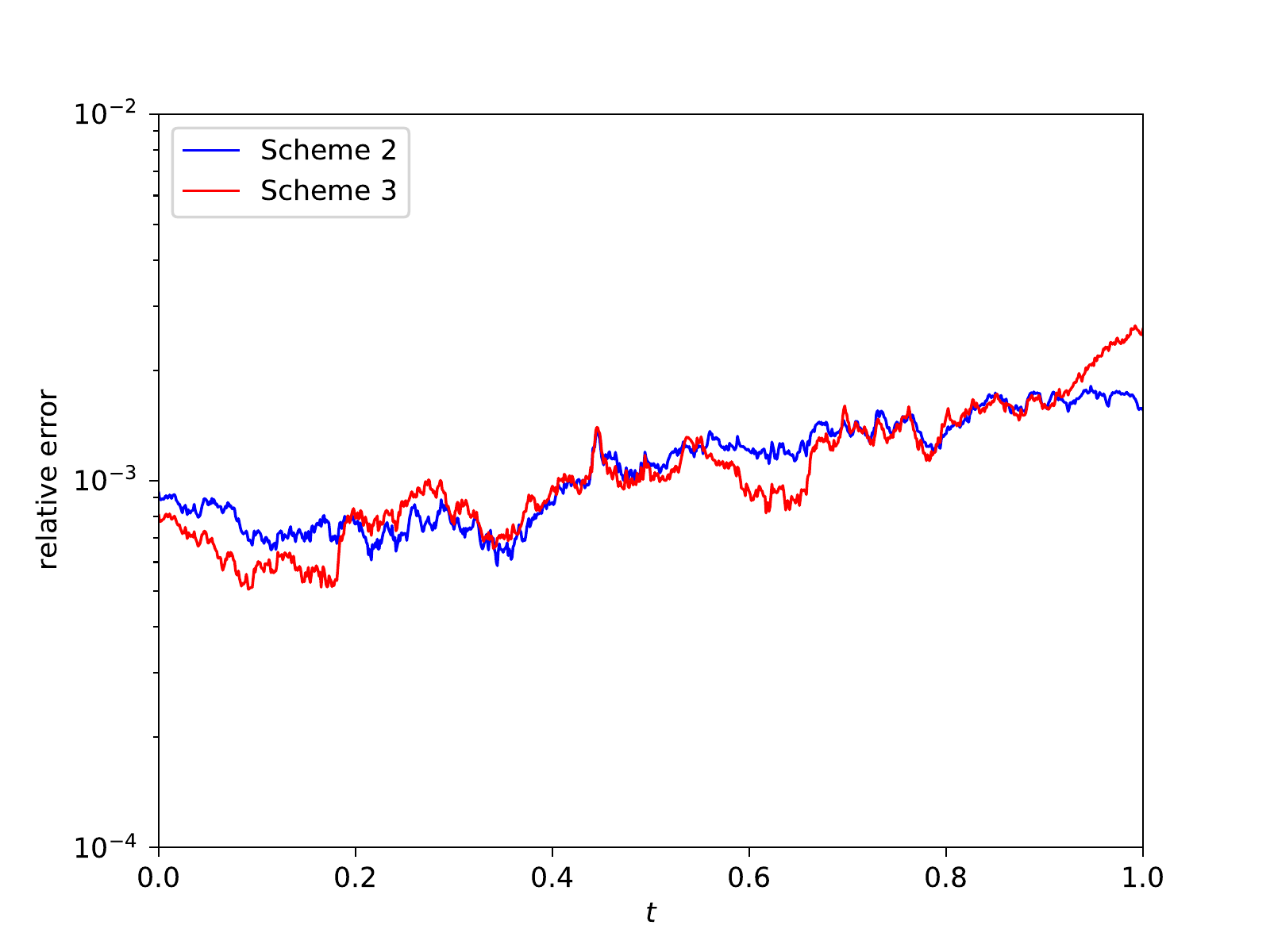}
\caption{Averaged relative error on $[0, T]$}
\end{subfigure}
\caption{Prediction of 8 test sample paths from training results of Scheme 2 and
Scheme 3, $N = 192$.}%
\label{plot-N192}%
\end{figure}

Fig. \ref{plot-N192} (a) (b) show the prediction of trained networks using
Scheme 2 and Scheme 3 with $N=192$ along 8 sampled test paths depicted in
Fig. \ref{plot-N192} (c), in comparison with the exact solution, where the
average error of the prediction is given in Fig. \ref{plot-N192} (d).

\subsubsection{Richardson extrapolation for higher order accuracy}
In Section \ref{numerical-S2S3} we have seen that Scheme 2 and Scheme 3 have the convergence behavior as the Euler--Maruyama scheme, so
we can assume that the truncation error may have the following asymptotic ansatz
\begin{equation}
u^{N}_\theta - u = C_1 N^{-\frac{1}{2}} + C_2 N^{-1} + O(N^{-\frac{3}{2}}),
\label{ansatz}
\end{equation}
where the leading term $C_1 N^{-\frac{1}{2}}$ dominates the error when $N$ is sufficiently large.
If this holds for both $u^{N}_\theta$ and $u^{4N}_\theta$ for some constants $C_1$ and $C_2$, then we can define an extrapolated solution
\begin{align}
u_{\mathrm{ex}}^{4N} {}&{}= 2 u_\theta^{4N} - u_\theta^{N} \nonumber \\
{}&{}= u -\frac{C_2}{2} N^{-1} + O(N^{-\frac{3}{2}})
\end{align}
as an improved approximation to the solution.

\begin{figure}[hbt]
\centering
\includegraphics[width=.45\linewidth]{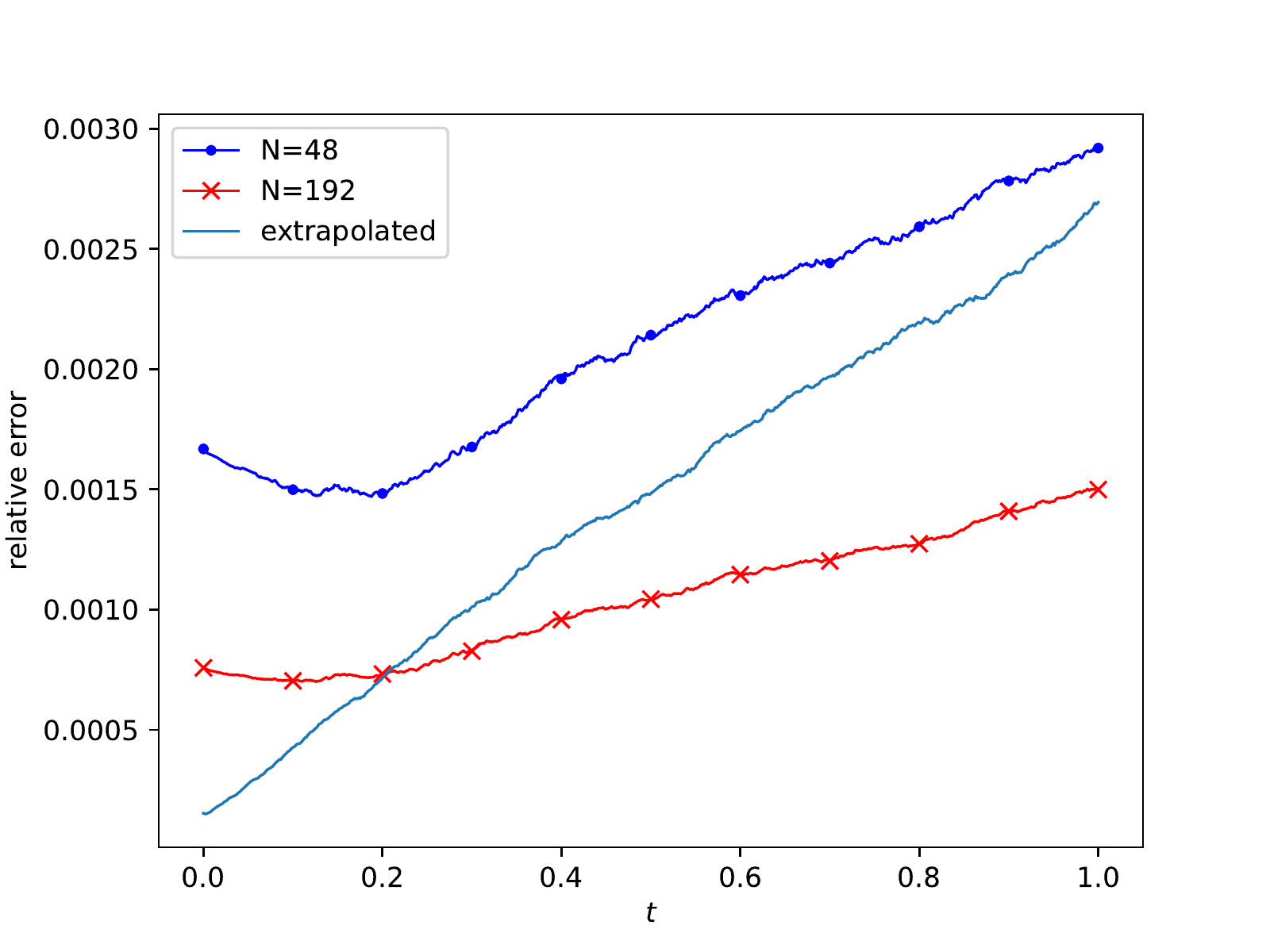}
\caption{Mean relative error of the extrapolation $u_\theta^{192}(t, X_t)$ for $0 \le t \le T$, using Scheme 2.}
\label{fig-extra}
\end{figure}

For the model problem (\ref{eq-BSB}), the Richardson extrapolation is valid for the approximation of $Y_0 = u(0, x_0)$ and $u(0,x)$ in a neighborhood near $x_0$, as shown by Table \ref{error-extra} and Fig. \ref{fig-extra}.
In terms of the accuracy of $Y_0$, by training the DNNs only with $N=12$ and $N=48$, the extrapolated result $u_{\mathrm{ex}}^{48} (0, x_0)$ has its accuracy outperforming those using $N=768$ which takes more than $10$ times longer time to train, when using both Scheme 2 and Scheme 3.
Due to training difficulties, the improvement for using extrapolation on $N=768$ is marginal, but still exists.

Note that the Richardson extrapolation approach usually may not work for the whole time interval along the entire sample paths. For instance, the values at $t=T$ are subject to explicit fitting of the terminal condition from the loss functions (\ref{scheme2loss}) and (\ref{loss-3}), so we cannot expect a general constant $C_1$ in (\ref{ansatz}) for $u^{N}_\theta(T, x)$ and $u^{4N}_\theta(T, x)$.
The result in Fig. \ref{fig-extra} shows that the extrapolation technique can be used for a time interval $0 \le t \le 0.1$.
\setlength{\heavyrulewidth}{1.5pt}
\setlength{\abovetopsep}{4pt}
\begin{table}[!htbp]
\centering
\begin{tabular}{*5c}
\toprule
{} &  \multicolumn{2}{c}{Scheme 2} & \multicolumn{2}{c}{Scheme 3}\\
\midrule
$N$   & $u^{N}_\theta$ & $u_{\mathrm{ex}}^{N}$ & $u^{N}_\theta$ & $u_{\mathrm{ex}}^{N}$\\
12   &  2.91e-03  &                     & 2.82e-03  &          \\
48   &  1.67e-03  & \textbf{4.29e-04}   & 1.13e-03  & 5.57e-04 \\
192  &  7.58e-04  & \textbf{1.53e-04}   & 8.43e-04  & 5.55e-04 \\
768  &  6.77e-04  & 5.97e-04            & 5.96e-04  & 3.49e-04 \\
\bottomrule
\end{tabular}
\caption{\label{error-extra}Relative error of $Y_0$ from the network approximation and extrapolation.}
\end{table}

\subsubsection{Region of validity of DNN $u_\theta(t,x)$ near $x_0$}
\begin{figure}[hbt!]
\centering
\begin{subfigure}{.45\textwidth}
\includegraphics[width=\linewidth]{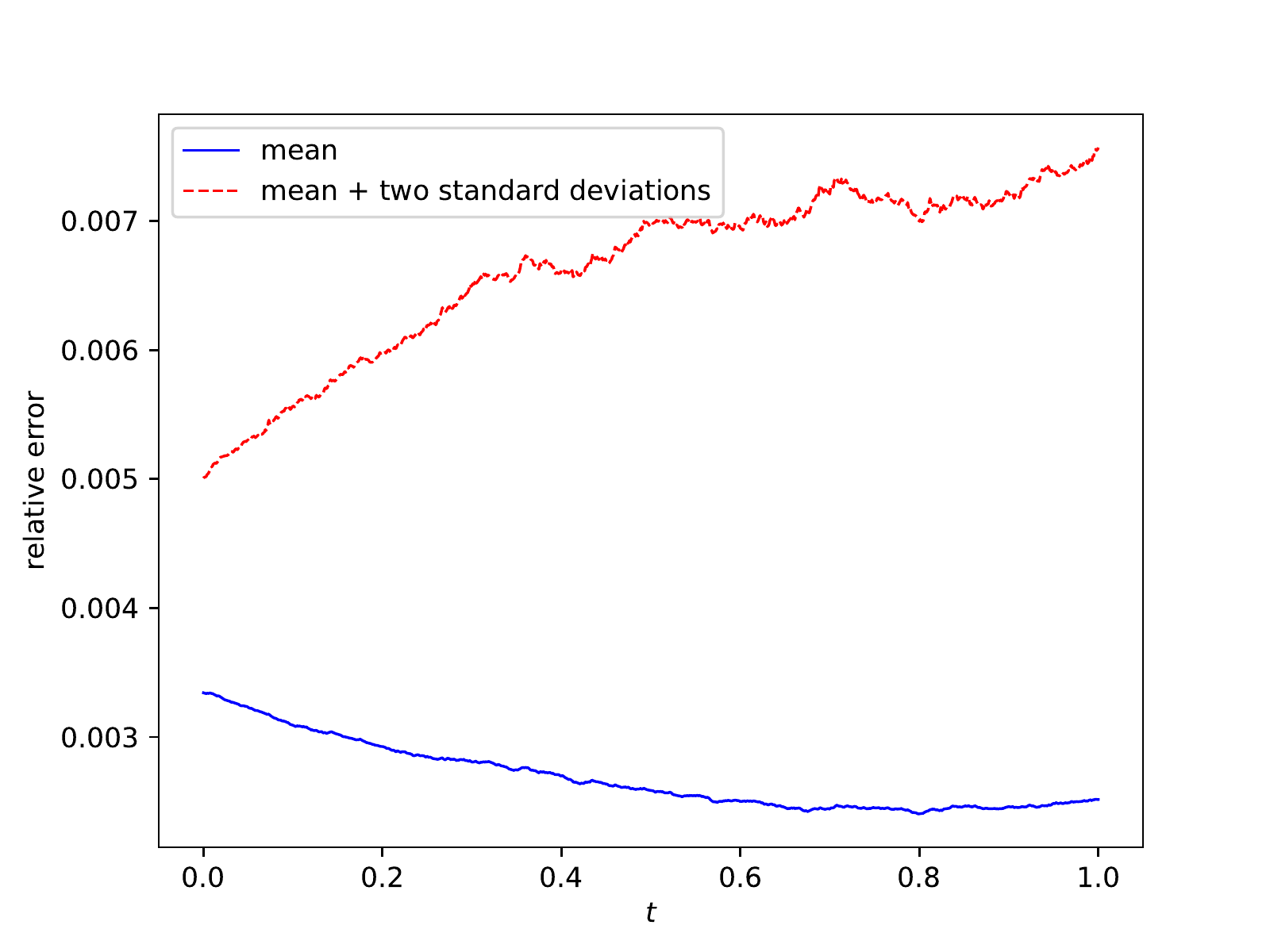}
\caption{Relative error, $R=0.25$}
\end{subfigure}
\begin{subfigure}{.45\textwidth}
\includegraphics[width=\linewidth]{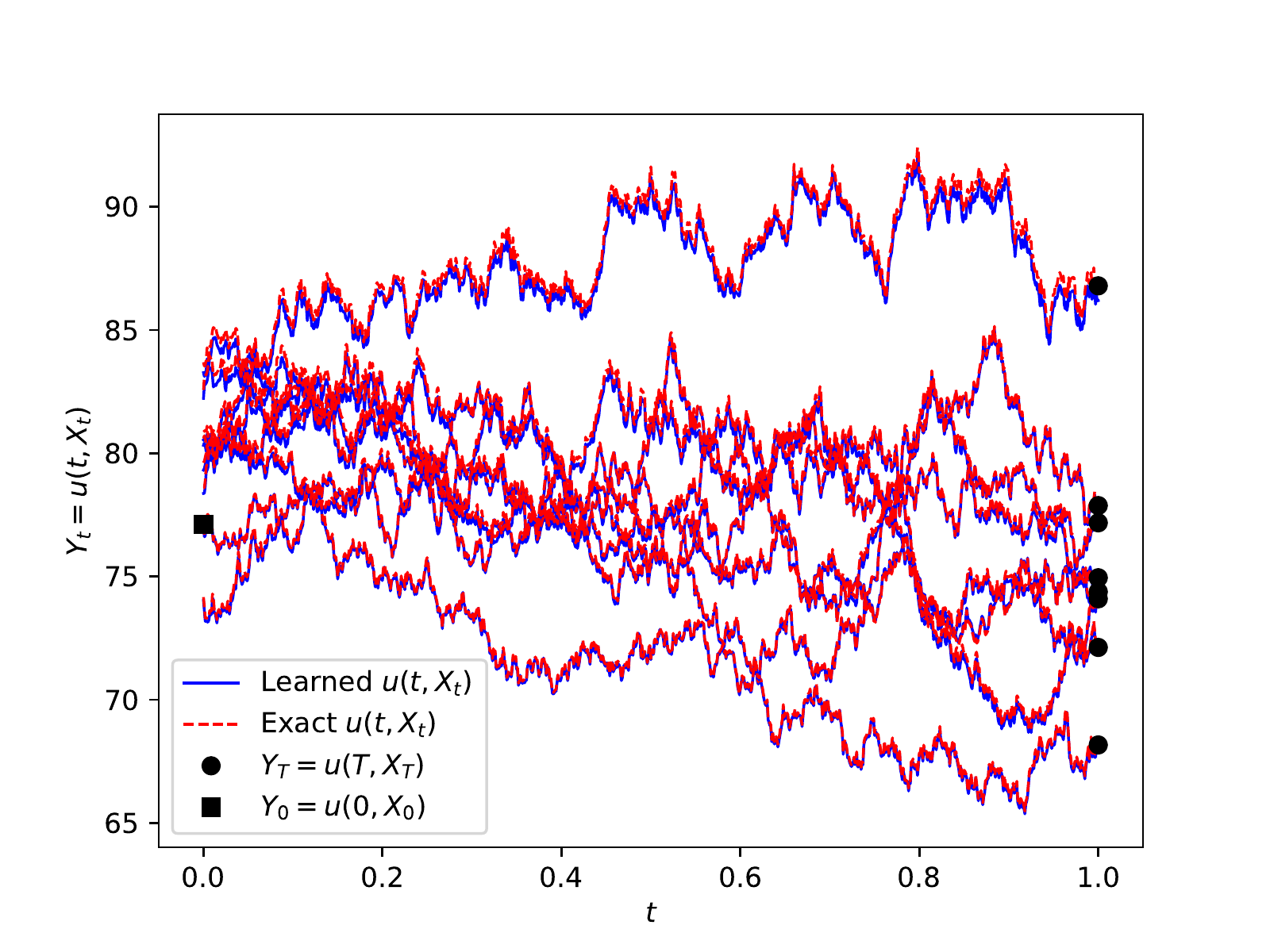}
\caption{Prediction of 8 sample paths}
\end{subfigure}
\hfill\begin{subfigure}{.45\textwidth}
\includegraphics[width=\linewidth]{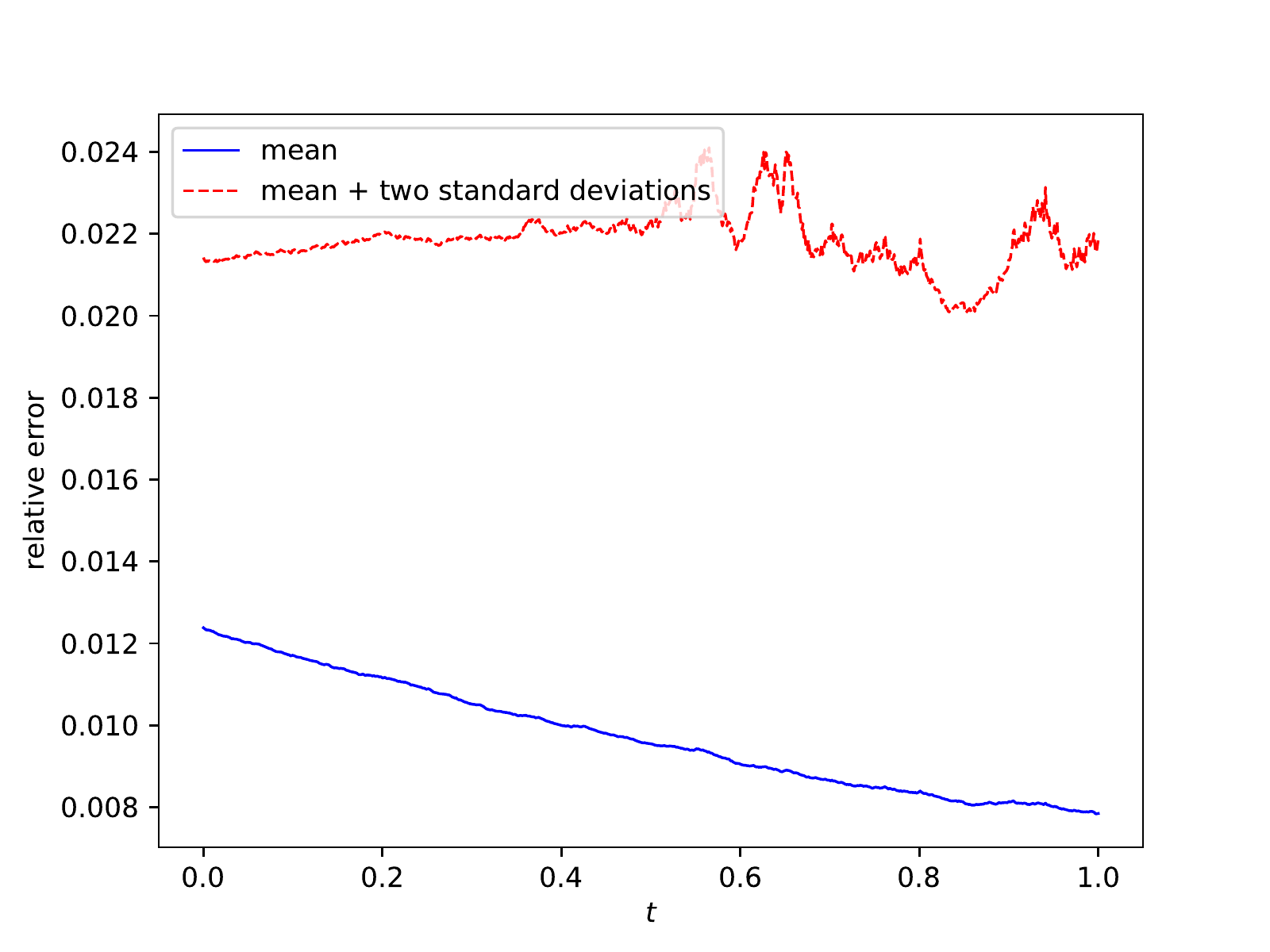}
\caption{Relative error, $R=0.5$}
\end{subfigure}
\begin{subfigure}{.45\textwidth}
\includegraphics[width=\linewidth]{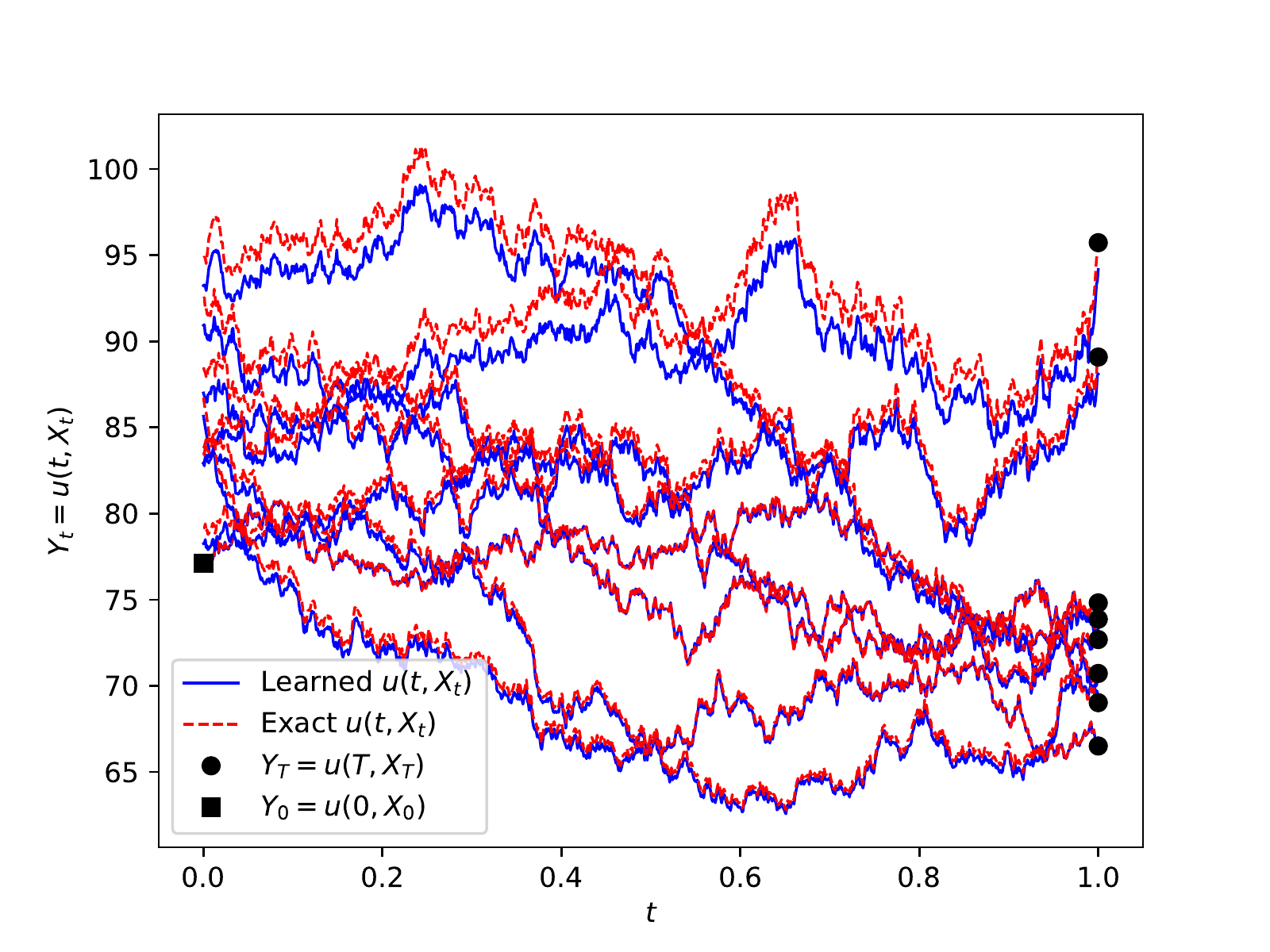}
\caption{Prediction of 8 sample paths}
\end{subfigure}
\caption{Training error verified with initial value $\tilde{x}_0$ from a neighborhood of $x_0$, using Scheme 2, $N = 192$.}%
\label{plot-nbhd}%
\end{figure}

In this section, we will verify the validity of the networks $u_\theta(t,x)$ in a region that are larger than the one sampled during the training process.
For this purpose, we randomly sample the initial value $X_0 = \tilde{x}_0$ from a \emph{cubic} neighborhood of $x_0$ with halved edge length $R$, i.e.,
\begin{equation}
(\tilde{x}_0)_j = (x_0)_j \cdot (1 + \varepsilon_j), \quad 1 \le j \le d=100,
\end{equation}
where $\varepsilon_j$ are i.i.d. random variables with uniform distribution on $(-R, R)$.
For the network trained with Scheme 2 and $N=192$, we compare the resulting error using the same measurement with $R=0.25$ and $R=0.5$, while keeping one sample starting exactly from $x_0$ (for the sake of plotting), see Fig. \ref{plot-nbhd}.
The averaged relative error is slightly larger at $t=0$ because during the training process these regions are less likely to be visited since we fixed the initial value for all training pathes at $X_0 = x_0$.
If we look at the overall maximum for $t \in [0, T]$, we can still have an averaged relative error of 0.34\% for $R=0.25$ and 1.25\% for $R=0.5$.
Also, it is noted that, in comparison with the non-perturbed result, the trained network fits the solution of the PDE better when $Y_t$ has a value below 80.

This result shows that the DNN we trained for $x=x_0$ is in fact can be used in a local neighbourhood around $x_0$ for the whole time interval $0 \le t \le T$.

\subsection{Multiscale DNN for the BSB equation with temporal oscillations}
\begin{figure}[ptbh]
\centering
\includegraphics[width=0.5\linewidth]{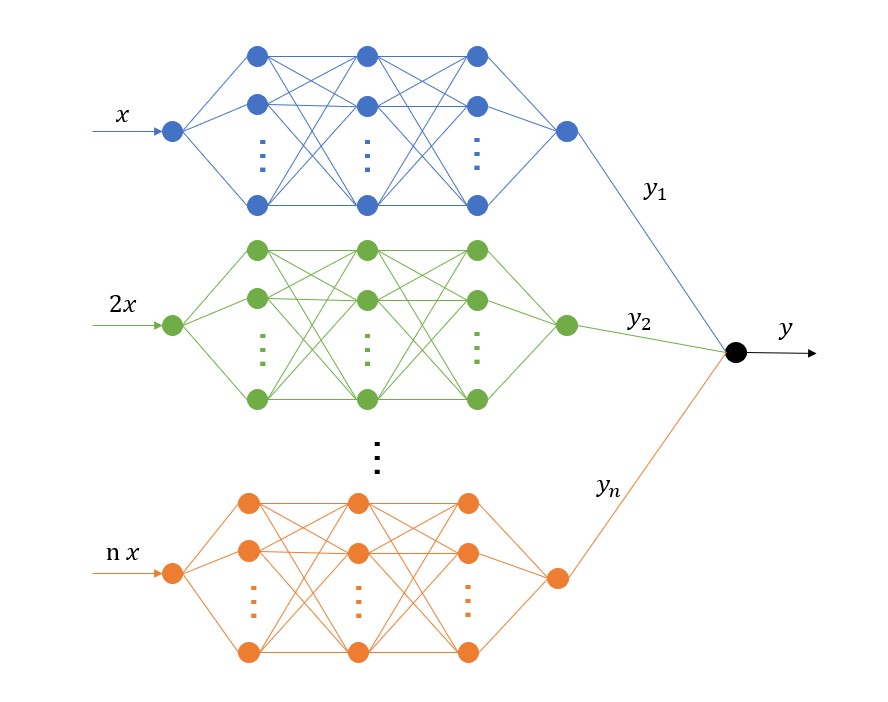}\caption{Illustration of a MscaleDNN.}%
\label{pic-net}%
\end{figure}

In a recent work \cite{msdnn}, a multi-scale DNN (MscaleDNN) was proposed, which consists of a series of parallel normal sub-networks, each of which receiving a scaled version of the input, and outputs of the sub-networks are combined to form the final output of the MscaleDNN (see Fig. \ref{pic-net}).
The individual sub-networks in the MscaleDNN with a scaled input is designed to approximate a segment of frequency content of the targeted function, and the effect of the scaling is to convert a specific high frequency content to a lower frequency range so that the learning can be accomplished more quickly, which is shown by the recent work \cite{msdnn} on the frequency dependence of the DNN convergence.

Fig. \ref{pic-net} shows the schematics of a MscaleDNN consisting of $n$ sub-networks.
Each scaled input passing through a fully-connected sub-network, which can be expressed in the formula (\ref{dnn}), here again we use the \emph{sine} function for the activation function, i.e.,
\begin{equation}
\sigma(x) = \sin(x).
\end{equation}
Mathematically, the final output of a MscaleDNN solution is represented by the following sum of sub-networks $f_{\theta^{n_i}}$ with network parameters denoted by $\theta^{n_i}$ (i.e. weight matrices and bias)
\begin{equation}
f(\bm{x}) \sim \sum_{i=1}^{M} \bm{W}_i^{[L]} f_{\theta^{n_i}} (\bm{\alpha}_i \cdot \bm{x}) + \bm{b}^{[L]},
\end{equation}
where $\bm{\alpha}_i$ is the chosen scale vector for the $i$-th sub-network in Fig. \ref{pic-net}.
For more details on the design of the MscaleDNN, refer to \cite{msdnn}.

For the input scales, the general idea is to adopt various scaling factors for different components of the input, depending on the complexity of the PDE to be solved.

\begin{figure}[htb]
\centering
\begin{subfigure}{.45\textwidth}
\includegraphics[width=\linewidth]{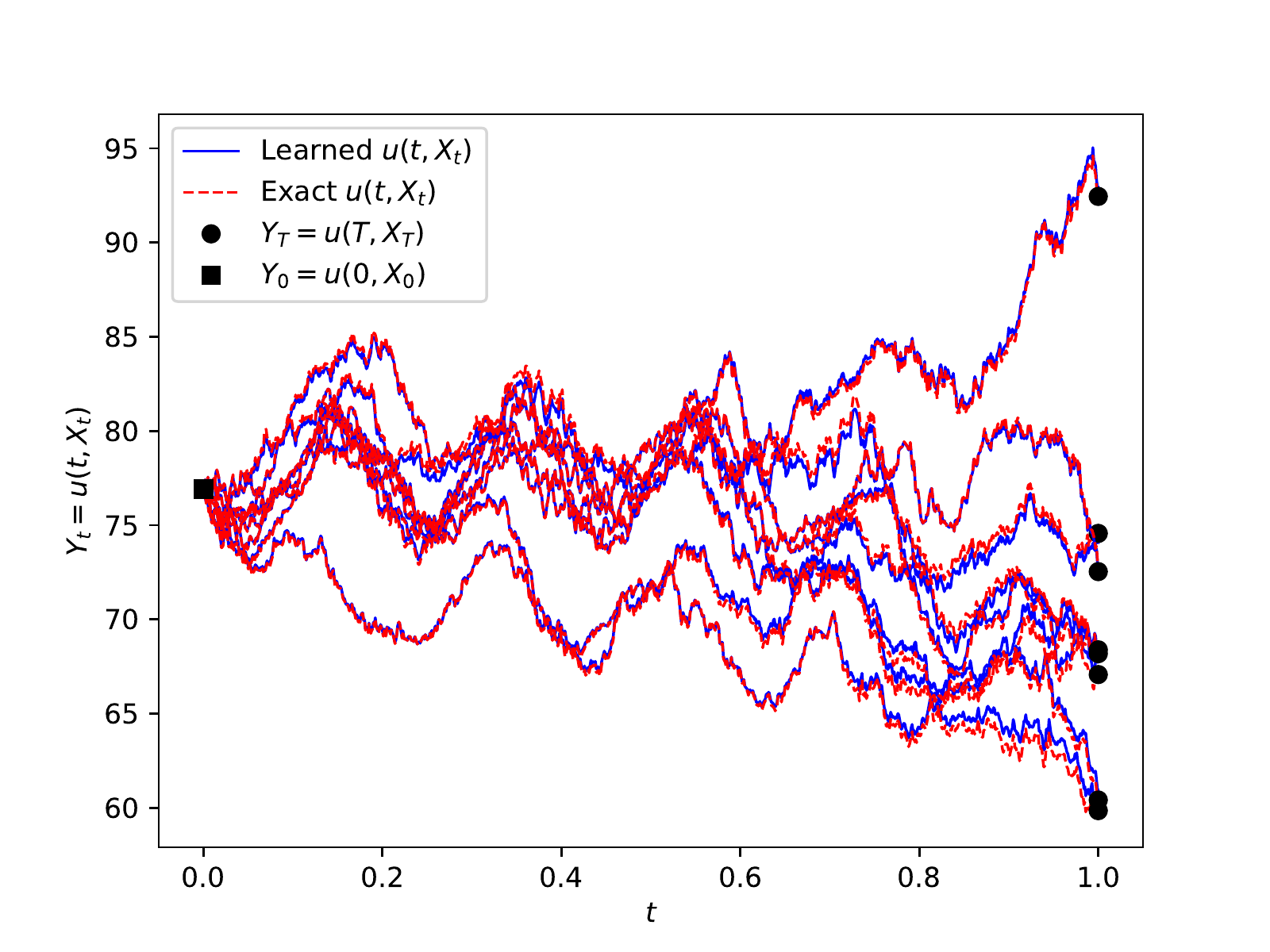}
\caption{Starting from $(0, x_0)$.}
\end{subfigure}
\begin{subfigure}{.45\textwidth}
\includegraphics[width=\linewidth]{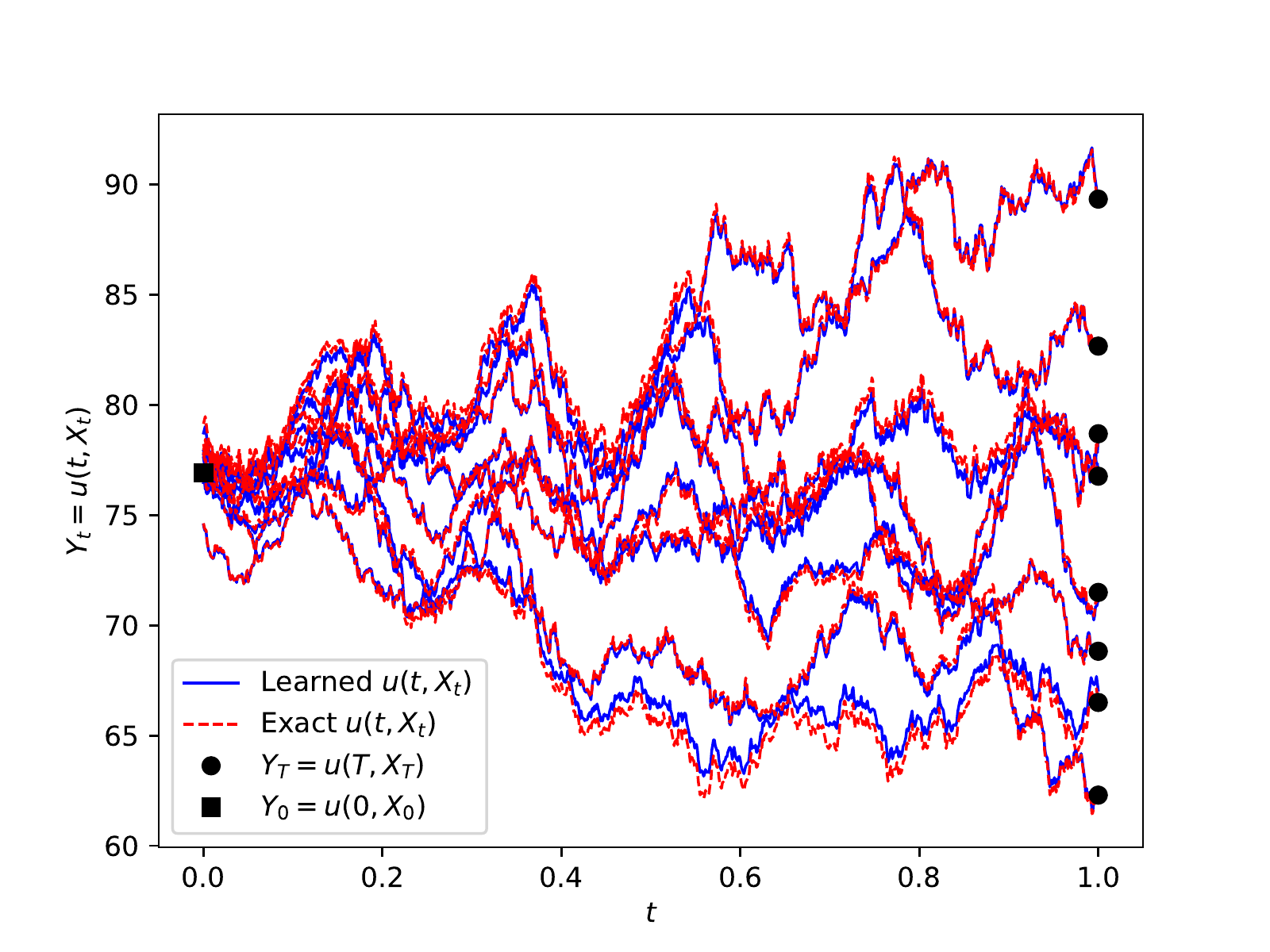}
\caption{Starting from a cubic neighborhood of $x_0$, $R=0.125$.}
\end{subfigure}
\caption{Prediction of 8 sample paths for problem with oscillation (\ref{eq-BSB-mod}), using the MscaleDNN with Scheme 2 and $N=192$.}
\label{fig-MS4-predict}
\end{figure}

The MscaleDNN is tested with the following model problem, modified from the BSB equation above with an oscillatory factor to effectively increase the training difficulty:
\begin{align}\label{eq-BSB-mod}
\begin{split}
\partial_t u + \frac{1}{2} \mathrm{Tr}[\sigma^2 \nabla \nabla u] &= \phi, \\
u(T, x) &= g(x),
\end{split}
\end{align}
where the dimension $d=100$, $T = 1.0$, $\sigma = 0.4$ and $r=0.05$ are unchanged parameters compared to (\ref{eq-BSB}),
\begin{equation}
g(x) = \lVert x \rVert^2 \left( 1 + \alpha \sin \left( \beta S_1 - \gamma T \right) \right),
\end{equation}
\begin{equation}
\phi(t, x, u, \nabla u) = r(u - \nabla u \cdot x) + \alpha e^{(r+\sigma^2)(T-t)} P(t, x),
\end{equation}
\begin{equation}
P(t, x) = (r \beta S_1 S_2 - \gamma S_2 + 2\sigma^2 \beta S_3) \cos(\beta S_1 - \gamma t) - \frac{\sigma^2 \beta^2}{2} S_2^2 \sin(\beta S_1 - \gamma t),
\end{equation}
where each $S_j = \sum_{i=1}^{d} x_i^j$, and $\alpha$, $\beta$ and $\gamma$ are parameters to be tuned.
The modified PDE (\ref{eq-BSB-mod}) has a solution
\begin{equation}
u(t, x) = e^{(r + \sigma^2) (T-t)} \lVert x \rVert^2 \left( 1 + \alpha \sin \left( \beta S_1 - \gamma t \right) \right),
\end{equation}
and corresponds to the FBSDEs
\begin{align}
\begin{split}
dX_{t}  &  =\sigma\mathrm{diag}(X_{t})dW_{t},\\
X_{0}  &  =x_0,\\
dY_{t}  &  =\left(r(Y_{t}-Z_{t}\cdot X_{t}) + \alpha e^{(r + \sigma^2) (T-t)} P(t, X_t) \right)dt+\sigma Z_{t}^{T}\mathrm{diag}%
(X_{t})dW_{t},\\
Y_{T}  &  =g(X_{T}).
\end{split}
\end{align}
We apply $\alpha = 0.025$, $\beta = 0.25$ and $\gamma = 32$ to the above equation.
During the training process, we use the same settings for the fully-connected DNN as in previous tests.
For the MscaleDNN, the network is divided into 4 sub-networks, each having 5 hidden layers with 64 neurons per layer, so that sizes of the networks in the comparison are matching.
The scaled inputs for the sub-networks are given by
\begin{equation}
(3^0 t,x), \quad (3^1 t, x), \quad (3^2 t, x), \quad (3^3 t, x),
\end{equation}
so that a wider range of frequency of $t$ can be captured with the MscaleDNN.
When applying Scheme 2 and $N=48$, the MscaleDNN halves the overall error compared to the fully-connected network.
\begin{figure}[htb!]
\centering
\begin{subfigure}{.45\textwidth}
\includegraphics[width=\linewidth]{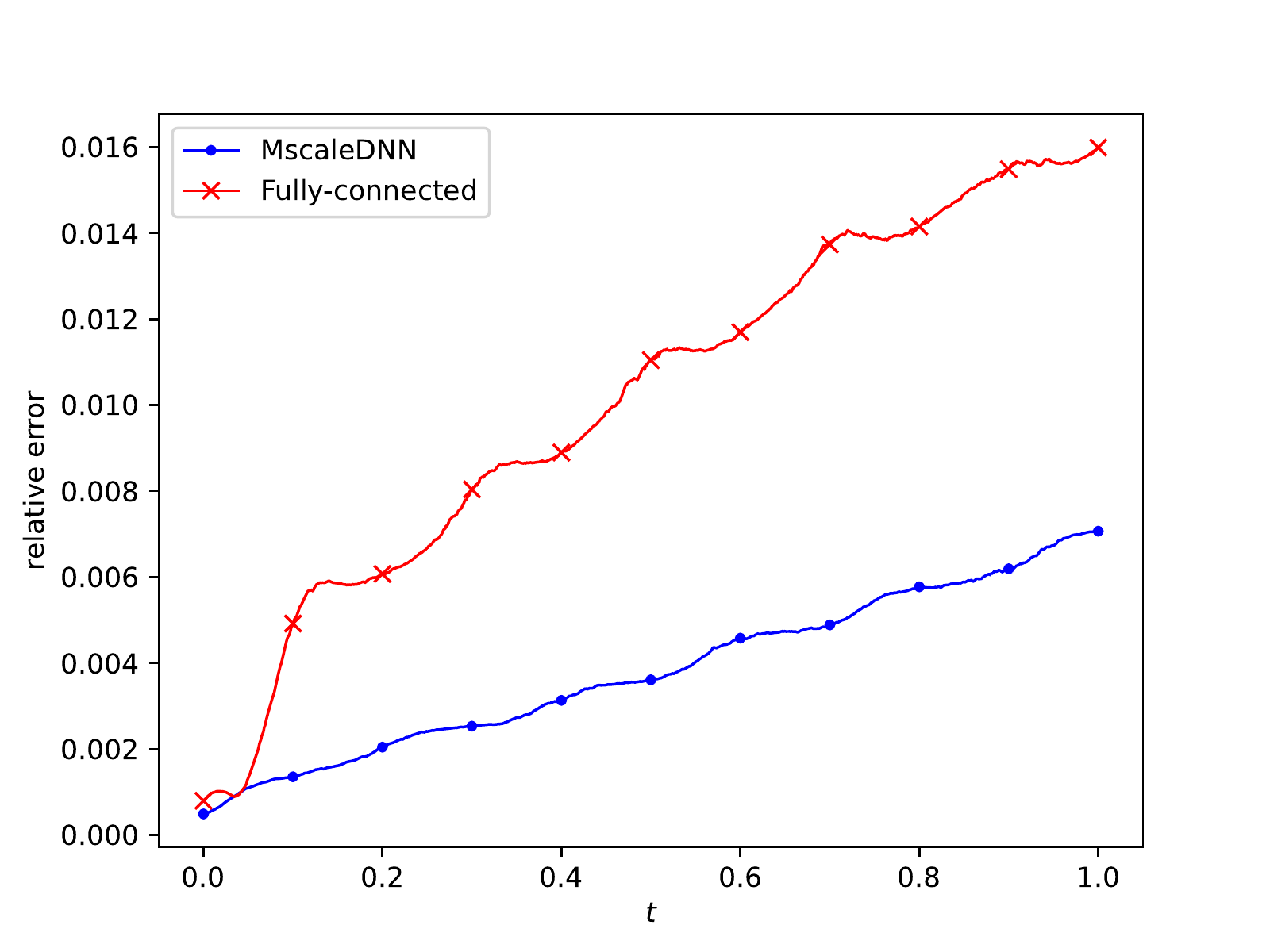}
\caption{Error Mean}
\end{subfigure}
\begin{subfigure}{.45\textwidth}
\includegraphics[width=\linewidth]{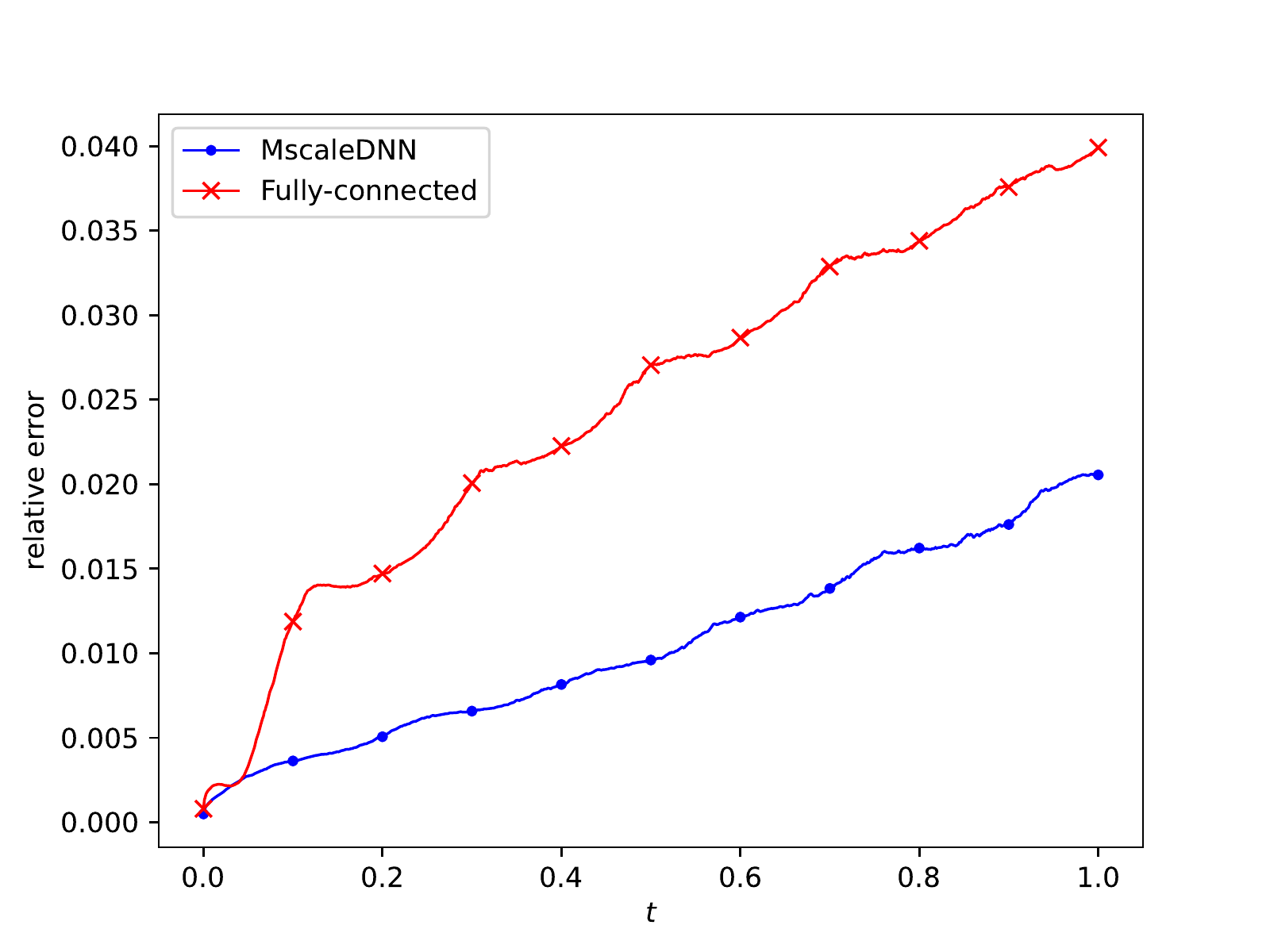}
\caption{Error Mean plus two SDs}
\end{subfigure}
\caption{Relative training error for fully-connected DNN and MscaleDNN for the model problem with oscillation, using Scheme 2 and $N=48$.}%
\label{fig-ms4}%
\end{figure}
One can also predict sample paths with better accuracy using the MscaleDNN, too, see Fig. \ref{fig-MS4-predict}.

\section{Conclusion}

In this paper, we have proposed two FBSDE based DNN algorithms for high
dimensional quasilinear parabolic equations. The key component of the
proposed algorithms is the loss function used, consisting of, in
addition to the terminal condition of the PDE, the pathwise difference of two
convergent stochastic processes from either the discretized SDEs or the PDEs
DNN solution. As the two stochastic processes converge to the same
stochastic process in the Pardoux--Peng theory, the new algorithms are able
to demonstrate nearly the half-order strong convergence of the underlying Euler--Maruyama scheme for the SDEs.
As a result,  the Richardson extrapolation method can be used, which confirms the convergence order of the DNN solutions and further enhances the resulting accuracy of the DNN solution of the PDE.
For PDEs with time oscillatory solutions, the MscaleDNN is shown to provide an enhancement of the resulting accuracy as well.

Future research will be done to improve the convergence of the networks and the overall schemes, including MscaleDNN for PDEs with spatially oscillatory solutions.


\end{document}